\documentclass[11pt,leqno]{article}

\usepackage{a4wide}
\usepackage{amssymb}
\usepackage[mathscr]{eucal}
\usepackage[all]{xy}

\newtheorem{theorem}{Theorem}
\newtheorem{lemma}[equation]{Lemma}
\newtheorem{proposition}[equation]{Proposition}
\newtheorem{corollary}[equation]{Corollary}
\newtheorem{conjecture}{Conjecture}
\newtheorem{definition}[equation]{Definition}
\newtheorem{remark}[equation]{Remark}

\newcommand{\resetL}{\setcounter{equation}{0}} 

\renewcommand{\AA}{\mathbb{A}}
\newcommand{\PP}{\mathbb{P}}
\newcommand{\CC}{\mathbb{C}}
\newcommand{\QQ}{\mathbb{Q}}
\newcommand{\ZZ}{\mathbb{Z}}
\newcommand{\RR}{\mathbb{R}}
\newcommand{\FF}{\mathbb{F}}
\newcommand{\Fq}{\mathbb{F}_q}
\newcommand{\Fp}{\mathbb{F}_p}
\newcommand{\NN}{\mathbb{N}}

\newcommand{\GG}{\mathbb{G}}

\newcommand{\HH}{\mathfrak{H}}
\newcommand{\fa}{\mathfrak{a}}

\newcommand{\fp}{\mathfrak{p}}

\newcommand{\fP}{\mathfrak{P}}
\newcommand{\fm}{\mathfrak{m}}

\newcommand{\fL}{\mathfrak{L}}

\newcommand{\OO}{\mathcal{O}}

\newcommand{\cP}{\mathcal{P}}
\newcommand{\II}{\mathcal{I}}
\newcommand{\JJ}{\mathcal{J}}
\newcommand{\TT}{\mathcal{T}}

\newcommand{\cD}{\mathscr{D}}

\newcommand{\End}{\mathrm{End}}

\newcommand{\Aut}{\mathrm{Aut}}

\newcommand{\HCM}{H_{\mathrm{CM}}}

\newcommand{\Gal}{\mathrm{Gal}}
\newcommand{\Pic}{\mathrm{Pic}}

\newcommand{\PGL}{\mathrm{PGL}_2}
\newcommand{\PSL}{\mathrm{PSL}_2}
\newcommand{\GL}{\mathrm{GL}_2}
\newcommand{\SL}{\mathrm{SL}_2}

\newcommand{\Stab}{\mathrm{Stab}}

\newcommand{\kinf}{k_\infty}

\newcommand{\bC}{{\bf C}}

\newcommand{\PGLA}{\mathrm{PGL}_2(A)}
\newcommand{\PGLk}{\mathrm{PGL}_2(k)}
\newcommand{\PGLki}{\mathrm{PGL}_2(k_\infty)}
\newcommand{\PSLA}{\mathrm{PSL}_2(A)}

\newcommand{\PSLki}{\mathrm{PSL}_2(k_\infty)}
\newcommand{\GLA}{\mathrm{GL}_2(A)}

\def\Endproof{{\hskip0pt\unskip\unskip\nobreak\hfil\penalty50
          \hskip1em\hbox{}\nobreak\hfil
           {$\square$}
          \parfillskip=0pt\finalhyphendemerits=0
          \par}\medskip}

\renewcommand{\matrix}[4]
{\left(\!\!\begin{array}{cc}#1\! &\! #2 \\ #3\! &\!
#4\end{array}\!\!\right)}

\begin{document}

\title{The Andr\'e-Oort conjecture for products of Drinfeld modular
curves}
\author{{\em Florian Breuer}, at Stellenbosch}
\date{\rule{30mm}{0.2pt}}
\maketitle

\paragraph{Abstract.}
Let $Z=X_1\times\cdots\times X_n$ be a product of Drinfeld modular
curves. We characterize those algebraic subvarieties $X \subset Z$
containing a Zariski-dense set of CM points, i.e. points
corresponding to $n$-tuples of Drinfeld modules with complex
multiplication (and suitable level structure). This is a
characteristic $p$ analogue of a special case of the Andr\'e-Oort
conjecture.

\section{Introduction}
\resetL

The aim of this paper is to prove an analogue of the Andr\'e-Oort
conjecture for the case of subvarieties of a product of Drinfeld
modular curves.

This conjecture states, in general,

\begin{conjecture}[Andr\'e-Oort]\label{AO}
Let $S$ be a Shimura variety, and $X\subset S$ an irreducible
algebraic subvariety. Then $X$ contains a Zariski-dense set of
special points if and only if $X$ is a subvariety of Hodge type.
\end{conjecture}
For the relevant definitions, and known results, we refer the
reader to \cite{Andre97,Belhaj01,Edixhoven98,Edixhoven01,
Edixhoven-Yafaev, Moonen95,Moonen98a,Moonen98b,
Oort97,Yafaev00,Yafaev01}.

The special case of Conjecture \ref{AO} of interest to us is the
following. Consider $\AA^n_{\CC}$ as the moduli space of
$n$-tuples of elliptic curves, by associating the tuple
$(E_1,\ldots,E_n)$ to the point $\big(j(E_1),\ldots,j(E_n)\big)$.
Then $\AA_{\CC}^n$ is a Shimura variety, and the special points
correspond to tuples for which each $E_i$ has complex
multiplication, hence they are also called {\em CM points}. The
subvarieties of Hodge type of $\AA_{\CC}^n$, which we call {\em
modular} subvarieties, are defined by imposing isogeny conditions
between some coordinates and setting other coordinates equal to CM
$j$-invariants. We make this definition precise below.

Denote by $\HH$ the Poincar\'e upper half-plane, on which the
group $\GL^+(\RR)$ acts. Then a point $\tau \in \HH$ corresponds
to the complex elliptic curve $E_{\tau} = \CC/\ZZ\tau\oplus\ZZ$.
Two points $\tau_1,\tau_2 \in \HH$ correspond to isogenous
elliptic curves if and only if $\tau_1 = \sigma(\tau_2)$ for some
$\sigma\in\GL^+(\QQ)$. We now let
$(\sigma_1,\ldots,\sigma_n)\in\GL^+(\QQ)^n$ and consider the map
\begin{eqnarray*}
\HH & \longrightarrow & \AA^n(\CC) \\
\tau & \longmapsto &
\Big(j(\sigma_1(\tau)),\ldots,j(\sigma_n(\tau))\Big).
\end{eqnarray*}
The image lies on an irreducible algebraic curve $Y\subset\AA^n$,
which we call a {\em modular curve} in $\AA^n$. Now we can define
the modular varieties in $\AA^n$ as products (up to permutation of
coordinates) of copies of $\AA^1$, CM points in $\AA^1$ and
modular curves in $\AA^m$ for $m\leq n$.

It is clear that a modular curve contains a dense (in the complex
topology) set of CM points (as all the coordinates are isogenous),
and hence so does a modular variety. Conjecture \ref{AO} claims
the converse, more precisely

\begin{conjecture}[Andr\'e-Oort for $\AA^n$]\label{AO2}
Let $X\subset\AA^n$ be an irreducible algebraic variety. Then $X$
contains a Zariski-dense set of CM points if and only if $X$ is a
modular variety in the above sense.
\end{conjecture}

Yves Andr\'e \cite{Andre98} has proved Conjecture \ref{AO2} for
$n=2$, and Bas Edixhoven \cite{Edixhoven98,EdixhovenPrep} has
shown that Conjecture \ref{AO2} holds for all $n$ if one assumes
the Generalized Riemann Hypothesis (GRH) for quadratic imaginary
fields.

The aim of this paper is to adapt Edixhoven's techniques to
function fields, and thereby prove Conjecture \ref{AO2} with
elliptic curves replaced by rank 2 Drinfeld modules.

More precisely, let $q$ be a power of the odd prime $p$ and set
$A=\Fq[T]$ and $k=\Fq(T)$. Denote by $\infty$ the place of $k$
with uniformizer $1/T$, let $\kinf = \Fq((1/T))$ be the completion
of $k$ at $\infty$, and let $\bC = \hat{\overline{k}}_{\infty}$
denote the completion of the algebraic closure of $\kinf$. (Here
$A$, $k$, $\kinf$ and $\bC$ play the roles of $\ZZ$, $\QQ$, $\RR$
and $\CC$, respectively). Then we may view $\AA_{\bC}^n$ as the
moduli space of $n$-tuples of rank 2 Drinfeld $A$-modules, via the
$j$-invariant, and a point $(x_1,\ldots,x_n)\in\AA^n(\bC)$ is
called a CM point if every $x_i$ is the $j$-invariant of a
Drinfeld module with complex multiplication. We prove the
following results.

\begin{theorem}\label{AOvarieties}
Assume that $q$ is odd. Let $X\subset \AA_{\bC}^n$ be an
irreducible variety. Then $X(\bC)$ contains a Zariski-dense subset
$S$ of CM points if and only if $X$ is a modular variety.
\end{theorem}

When $X$ is a curve, we have an effective result.

\begin{theorem}\label{AOcurves}
Assume that $q$ is odd. Let $d,m$ and $n$ be given positive
integers, and $g$ a given non-negative integer. Then there exists
an effectively computable constant $B=B(n,m,d,g)$ such that the
following holds. Let $X$ be an irreducible algebraic curve in
$\AA_{\bC}^n$ of degree $d$, defined over a finite extension $F$
of $k$ of degree $[F:k] = m$ and genus $g(F)=g$. Then $X$ is a
modular curve if and only if $X(\bC)$ contains a CM point of
height at least $B$.
\end{theorem}

As level structures play no role here, one may replace $\AA^n =
(\AA^1)^n$ by a product $X_1\times\cdots\times X_n$ of Drinfeld
modular curves, and obtain a similar result (see Corollaries
\ref{toAOcurves} and \ref{toAOvarieties} for the exact
statements). The definition of modular curves and modular
varieties, in $\AA^n$ or in $X_1\times\cdots\times X_n$, will be
given in \S\ref{ModularDefSection}.

The proofs of Theorems \ref{AOvarieties} and \ref{AOcurves} may be
divided into two parts. Firstly, in the topological part
(\S\ref{Hecke}) one shows that a variety which is stabilized by a
certain Hecke operator must be modular. In this part we follow an
approach similar to Edixhoven's, but the translation from number
fields to function fields is not automatic, as problems specific
to finite characteristic arise (e.g. one can no longer use
arguments from Lie theory, $\kinf$ has many non-trivial
automorphisms, and $\bC$ has infinite dimension over $\kinf$).

In the second, arithmetic part (\S\ref{Heights}), one shows that
varieties containing suitable CM points are stabilized by certain
Hecke operators, and the translation is easier. Here one uses GRH
twice, once for a strong version of the \v{C}ebotarev Theorem, and
once to obtain effective bounds on class numbers of quadratic
imaginary fields (one only needs GRH for such bounds if one wants
an {\em effective} result, such as Theorem \ref{AOcurves}, the
classical analogue of which appears in \cite{Breuer01}). But as
we're working over function fields, GRH is known (Hasse-Weil), so
our results are unconditional.

We have assumed that $q$ is odd for technical reasons (notably
concerning the arithmetic of quadratic extensions of $k$), but we
expect a similar result to hold in characteristic $2$.

In the rest of this introduction we will briefly recall some facts
about Drinfeld modules and Drinfeld modular curves, to fix
notation, and we will define the notion of modular varieties in
\S\ref{ModularDefSection}.

\paragraph{Acknowledgments.} The results presented here made up my Ph.D thesis
at l'Universit\'e Denis Diderot (Paris 7), and I am
deeply indebted to my supervisor, Marc Hindry, for his cheerful
advice and guidance. I would also like to thank Bas Edixhoven for
his many patient explanations, and for making a preliminary
version of \cite{EdixhovenPrep} available to me. The idea of
replacing elliptic curves by Drinfeld modules in the Andr\'e-Oort
conjecture was first suggested to me by Hans-Georg R\"uck, and I
am also grateful to Henning Stichtenoth for providing me with
Proposition \ref{Boundh}, which allowed me to remove the condition
$q\geq 5$. Lastly, I wish to thank the National Center for Theoretical
Sciences in Hsinchu, Taiwan, and the Max-Planck-Institut f\"ur 
Mathematik in Bonn, Germany, for their hospitality.

\subsection{Drinfeld modules}

We give here a very brief introduction to Drinfeld modules, the
aim being rather to fix our notation than to initiate the reader
in this fascinating topic. For details, we refer the reader to
\cite{Goss96} and \cite{Hayes92}.

Let $\tau$ denote the $q$-th power Frobenius acting on the
additive group $\GG_{a,\bC}$. Then the $\Fq$-linear endomorphisms
of $\GG_{a,\bC}$ are given by
$\End_{\Fq}(\GG_{a,\bC}) = \bC\{\tau\},$
the ring of {\em twisted polynomials}, i.e. non-commutative
polynomials in $\tau$ subject to the relations $\tau x = x^q \tau$
for all $x\in\bC$. Then a {\em Drinfeld $A$-module of rank $r$}
(and defined over $\bC$) is an injective ring homomorphism
\begin{eqnarray*}
\phi : A & \longrightarrow & \bC\{\tau\} \\
       a & \longmapsto     & \phi_a = a\tau^0 + a_1\tau + \ldots + a_n\tau^n,\; a_n\neq 0,\; n=r\deg(a).
\end{eqnarray*}

A morphism of Drinfeld modules, written\footnote{Beware, this is
not a function between sets!} $f : \phi \rightarrow \phi'$, is an
element $f\in\bC\{\tau\}$ such that $f\phi_a = \phi'_af$ for all
$a\in A$. If $f$ is non-zero then $\phi$ and $\phi'$ have the same
rank, and we call $f$ an isogeny. $f$ is an isomorphism if and
only if $f\in\bC^*$. The isogeny $f$ is called {\em cyclic of
degree} $N\in A$ if $\ker(f) \cong A/NA$ as $A$-modules.

There is also an analytic construction of Drinfeld modules,
similar to the construction of complex elliptic curves as
quotients of $\CC$ by a lattice. A {\em lattice} of rank $r$ in
$\bC$ is a discrete $A$-submodule $\Lambda$ of $\bC$ such that
$\kinf\Lambda$ has dimension $r$ over $\kinf$. As $\bC$ has
infinite dimension over $\kinf$, it contains lattices of any rank
(unlike $\CC$). Then there is an equivalence between the
categories (Drinfeld modules of rank $r$ over $\bC$, morphisms)
and (Lattices of rank $r$ in $\bC$, homotheties).

From now on, by a Drinfeld module $\phi$ we will always mean a
Drinfeld $A$-module of rank $r=2$.

Clearly $\End(\phi) = \{ f\in\bC\{\tau\} \;|\; f\phi_a = \phi_af,
\forall a\in A\}$ is the centralizer of $\phi(A)$ in
$\bC\{\tau\}$. Generically, $\End(\phi)\cong A$, but sometimes
$\End(\phi)$ is strictly larger than $A$, and we say that $\phi$
has {\em complex multiplication} (CM). In this case, the
endomorphism ring is of the form $\OO=A[\sqrt{d}]$, for some
non-square $d\in A$. Write $d=f^2D$, with square-free $D\in A$,
then $\OO$ is an order of conductor $f$ in the {\em CM field}
$K=k(\sqrt{D})$. $K/k$ is an {\em imaginary} quadratic extension,
which means that the place $\infty$ does not split in $K/k$.
Equivalently, $K$ has no embedding into $\kinf$, hence the
terminology. We distinguish two cases. Either $\infty$ ramifies in
$K/k$, in which case $\deg(D)$ is odd, or $\infty$ is inert in
$K/k$, in which case $\deg{D}$ is even and the leading coefficient
of $D$ is not a square in $\Fq$.

A Drinfeld module $\phi$ is uniquely determined by $\phi_T =
T\tau^0 + g\tau + \Delta\tau^2$, where $g,\Delta\in\bC$ and
$\Delta\neq 0$. We define the $j$-invariant of $\phi$ by
$j=j(\phi)=g^{q+1}/\Delta$, and one verifies easily that two
Drinfeld modules $\phi$ and $\phi'$ are isomorphic if and only if
$j(\phi)=j(\phi')$. Moreover, $\phi$ is isomorphic to a Drinfeld
module defined over $k(j)$.

The $j$-invariant induces a bijection between the set of
isomorphism classes of Drinfeld modules and $\AA^1(\bC)$, hence we
view $\AA_{\bC}^n$ as the moduli space of $n$-tuples of Drinfeld
modules. A point $x=(x_1,\ldots,x_n)\in\AA^n(\bC)$ is called a
{\em CM point} if each $x_i$ is the $j$-invariant of a CM Drinfeld
module.
As in the classical case, CM points have remarkable arithmetical
properties (see \cite{Gekeler82}):

\begin{theorem}[CM Theory]\label{CMtheory}
Let $\phi$ be a CM Drinfeld module, with $\OO=\End(\phi)$ an order
of conductor $f$ in the CM field $K=k(\sqrt{D})$. Then
$j=j(\phi)\in k^{sep}$ is integral over $A$, and $K(j)$ is the
ring class field of $\OO$. This means that $K(j)/K$ is unramified
outside $f$, split completely at the unique place $\infty$ of $K$,
and $\Gal(K(j)/K) \cong \Pic(\OO)$ via class-field theory.
\end{theorem}

In particular, suppose that $\fp\in A$ is a prime which splits in
$K/k$ and does not divide $f$, we say that $\fp$ splits in $\OO$.
Let $\sigma_{\fp} = (\fp, K(j)/K)$ be the associated Frobenius
element. Then there is a cyclic isogeny of degree $\fp$ from
$\phi$ to $\phi^{\sigma_{\fp}}$.

\subsection{Drinfeld modular curves}

Let $\Omega=\bC \smallsetminus \kinf$ denote Drinfeld's upper
half-plane, which plays the role of $\HH^{\pm}$ in the classical
case. The group $\PGLki$ acts on $\Omega$ by fractional linear
transformations, but unlike the classical case, this action is not
transitive, as $\bC$ has infinite dimension over $\kinf$. A point
$z\in\Omega$ is called {\em quadratic} if $[\kinf(z) : \kinf]=2$,
in which case the stabilizer of $z$ in $\PGLki$ is a one
dimensional Lie group over $\kinf$. Otherwise we call $z$ {\em
non-quadratic}, and its stabilizer is trivial.

Similarly to the classical case, the quotients of $\Omega$ by
congruence subgroups of $\PGLA$ give rise to affine Drinfeld
modular curves, which may be compactified by adding finitely many
cusps. See \cite{GekelerDMC} for details.

A point $z\in\Omega$ gives rise to a Drinfeld module $\phi^z$
associated to the lattice $\langle 1,z \rangle$. Notice that
$\phi^z$ has CM if and only if $[k(z):k]=2$. The $j$-invariant
induces a rigid analytic isomorphism $j:\PGLA\diagdown\Omega
\stackrel{\sim}{\rightarrow} \AA^1(\bC)$.

We will here define three special Drinfeld modular curves.
Let $N\in A$ and set
\begin{eqnarray*}
\Gamma(N) & = & \{ \gamma\in\GL(A) \;|\; \gamma \equiv 1 \bmod N\}/Z(\Fq^*) \subset\PGLA, \\
\Gamma_0(N) & = & \left\{ \matrix{a}{b}{c}{d} \in \GL(A) \;|\; c
\equiv 0 \bmod N \right\}/Z(\Fq^*)
      \subset\PGLA, \\
\Gamma_2(N) & = & \{ \gamma\in\GL(A) \;|\; (\gamma \bmod N)\in
Z\big((A/NA)^*\big)\}/Z(\Fq^*) \subset\PGLA.
\end{eqnarray*}
Here $Z(R^*)\cong R^*$ denotes the subgroup of scalar matrices in
$\GL(R)$, for a ring $R$. Then we define the curves\footnote{There
also exist Drinfeld modular curves denoted by $Y_1(N)$, which is
why we introduce the notation $Y_2(N)$.}
\[
Y(N) = \Gamma(N) \diagdown\Omega, \quad Y_0(N) = \Gamma_0(N)
\diagdown\Omega, \quad Y_2(N) = \Gamma_2(N) \diagdown\Omega.
\]

The curve $Y_0(N)$ is the coarse moduli space parameterizing
isomorphism classes of triples $(\phi, \phi', f)$, where $\phi$
and $\phi'$ are Drinfeld modules and $f : \phi \rightarrow \phi'$
is a cyclic isogeny of degree $N$. We have $Y(1) = Y_0(1) = Y_2(1)
= \AA^1$.

\begin{proposition}\label{GaloisCover}
The curve $Y_2(N)$ covers $Y_0(N)$ and is Galois over $Y(1)$.
Suppose that $N$ is square-free and that every prime factor of $N$
has even degree. Then $\Gal\big(Y_2(N)/Y(1)\big) \cong
\PSL(A/NA)$.
\end{proposition}
Notice that we use the definition
\[
\PSL(R):=\SL(R)/(Z(R^*)\cap\SL(R)) \cong \SL(R)/\{x\in R^* \;|\;
x^2=1\}
\]
for a ring $R$.

\paragraph{Proof.}Firstly, it is clear that $Y_2(N)$ covers $Y_0(N)$, as $\Gamma_2(N)\subset\Gamma_0(N)$.

We note that $Y(N)/Y(1)$ is Galois \cite{Gekeler84} with Galois
group
\[
\Gal\big(Y(N)/Y(1)\big)  \cong  G(N)/Z(\Fq^*),
\]
where we have set
\[
G(N) = \{\alpha\in\GL(A/NA) \;|\; \det(\alpha) \in \Fq^* \}.
\]

Next, we consider the coverings $Y(N) \rightarrow Y_2(N)
\rightarrow Y(1)$. Here $Y_2(N)$ corresponds to the normal
subgroup $H(N)/Z(\Fq^*)$ of $\Gal\big(Y(N)/Y(1)\big)$, where
\[
H(N) = Z\big((A/NA)^*\big)\cap G(N).
\]
Hence $Y_2(N)$ is Galois over $Y(1)$, with Galois group
\[
\Gal\big( Y_2(N)/Y(1) \big)  \cong G(N)/H(N) \subset \PGL(A/NA).
\]
This is the subgroup of $\PGL(A/NA)$ of those elements with
determinant in $\Fq^*$. As $N$ is square-free, composed of prime
factors of even degree, it follows that every $\alpha\in\Fq^*$ is
a square in $A/NA$. Hence $G(N)/H(N) = \PSL(A/NA)$. \Endproof

\subsection{Pure modular curves in $\AA^n$}\label{CurveSection}

For a subset $I\subset\{1,\ldots,n\}$ we denote by $p_I : \AA^n
\rightarrow \AA^I$ the projection onto the coordinates listed in
$I$. We also write $p_{i,j}=p_{\{i,j\}}$.

We map $Y_0(N)$ into $\AA^2$ by sending the pair $(\phi,\phi')$ of
isogenous Drinfeld modules to the point $(j(\phi),j(\phi'))$. The
image, which we denote by $Y'_0(N)$, is the locus of an
irreducible polynomial $\Phi_N(t_1,t_2) \in A[t_1,t_2]$ (see
\cite{Bae92}). This polynomial is symmetrical and of degree
$\psi(N) = |N|\prod_{\fp|N}(1+|\fp|^{-1})$ in $t_1$ and $t_2$.

Two points $z_1,z_2 \in \Omega$ correspond to isogenous Drinfeld
modules $\phi^{z_1}$ and $\phi^{z_2}$ if and only if $z_1 =
\sigma(z_2)$ for some $\sigma\in\PGLk$. In this case, $\phi^{z_1}$
and $\phi^{z_2}$ are linked by a cyclic isogeny of degree
$N=\det(m\sigma)$, where $m\in A$ is chosen in such a way that the
entries of $m\sigma$ are in $A$ and have no factor in common. We
call this $N$ the {\em degree} of $\sigma$ (it is unique up to the
square of an element of $\Fq^*$).

Let $(\sigma_1,\ldots,\sigma_n)\in\PGLk^n$ and consider the map
\begin{eqnarray}
\rho : \Omega & \longrightarrow & \AA^n(\bC) \label{immersion}\\
z & \longmapsto &
\big(j(\sigma_1(z)),\ldots,j(\sigma_n(z))\big).\nonumber
\end{eqnarray}
The image lies on an irreducible algebraic curve $Y\subset\AA^n$,
which we call a {\em pure modular curve}. When $n=2$, we obtain
again the curve $Y'_0(N)$, where $N$ is the degree of
$\sigma_2\sigma_1^{-1}$. We note that a curve $Y\subset\AA^n$ is a
pure modular curve if and only if $p_{i,j}(Y) = Y'_0(N_{i,j})$ for
some $N_{i,j}\in A$ and every pair of coordinates $i\neq j$.

Clearly $(\sigma_1,\ldots,\sigma_n)$ and
$(\sigma'_1,\ldots,\sigma'_n)$ in $\PGLk^n$ define the same
modular curve via the map (\ref{immersion}) if and only if there
exist $\sigma\in\PGLk$ and $(\gamma_1,\ldots,\gamma_n) \in\PGLA^n$
such that $\gamma_i\sigma_i\sigma = \sigma'_i$ for $i=1,\ldots,n$.
So the set of pure modular curves in $\AA^n$ is in bijection with
the double cosets
\begin{equation}\label{doublecoset1}
\big(\PGL(A)\diagdown\PGL(k)\big)^n\diagup\PGL(k),
\end{equation}
where the $\PGLk$ acts diagonally on $\PGLk^n$. We may transpose
the actions, and write (\ref{doublecoset1}) as the restricted
product, over all primes $\fp\in A$, of
\[
\PGL(k_{\fp})\diagdown \big(\PGL(k_{\fp}) \diagup
\PGL(A_{\fp})\big)^n,
\]
where $k_{\fp}$ and $A_{\fp}$ denote the completions at $\fp$ of
$k$ and $A$, respectively. Now each \linebreak $\PGL(k_{\fp})/\PGL(A_{\fp})$
may be identified with the Br\^uhat-Tits tree $\TT_{\fp}$ of
$\PGL(k_{\fp})$ (see \cite{SerreTrees}). We recall that the
vertices of $\TT_{\fp}$ correspond to homothety classes of
$A_{\fp}$-lattices in the vector space $k_{\fp}^2$, and two
vertices $v_1$ and $v_2$ are adjacent if we may find
representative lattices $L_i\in v_i$ such that $L_1/L_2 \cong
A_{\fp}/\fp A_{\fp}$. So a pure modular curve $Y$ corresponds to
an $n$-tuple $(v_{\fp,1},\ldots,v_{\fp,n})$ of vertices of
$\TT_{\fp}$, up to $\PGL(k_{\fp})$-action, for every prime $\fp\in
A$, with the condition that the vertices
$v_{\fp,1},\ldots,v_{\fp,n}$ coincide for almost all $\fp$.

We now give a special description of the pure modular curves in
$\AA^3$.
An {\em end} of $\TT_{\fp}$ is an equivalence class of infinite paths 
$\xymatrix@C=10pt@M=-1pt{\bullet\ar@{-}[r]&\bullet\ar@{-}[r]&\bullet\ar@{-}[r]&\;\cdots}$
of distinct vertices in $\TT_ {\fp}$, two such paths being
equivalent if they differ by a finite subgraph. The set of ends of
$\TT_ {\fp}$ is in a natural bijection with $\PP^1(k_{\fp})$.
Given a triple of distinct ends $(E_1,E_2,E_3)$ we may define a
unique vertex $v_c$ of $\TT_{\fp}$, called the {\em center} of
$(E_1,E_2,E_3)$, such that each $E_i$ is represented by a path
starting at $v_c$. Similarly, for a triple of (not necessarily
distinct) vertices $(v_1,v_2,v_3)$, we may define the center $v_c$
as the unique vertex such that there exist disjoint paths from
$v_c$ to each $v_i$. As the action of $\PGL(k_{\fp})$ on
$\PP^1(k_{\fp})$ is $3$-transitive, we may map any given triple
$(E_1,E_2,E_3)$ of distinct ends with center $v_c$ to any other
given triple $(E'_1,E'_2,E'_3)$ with center $v'_c$ via an element
of $\PGL(k_{\fp})$. It follows that the $\PGL(k_{\fp})$-class of
triples of vertices $(v_1,v_2,v_3)$ of $\TT_{\fp}$ is uniquely
determined by the triple $(n_1,n_2,n_3)$ of distances from the
vertices $v_i$ to the center of $(v_1,v_2,v_3)$.

Let $Y \subset \AA^3$ be a pure modular curve, corresponding to
triples of vertices $(v_{\fp,1},v_{\fp,2},v_{\fp,3})$ in $\TT_
{\fp}$ for each prime $\fp\in A$. To each triple we associate the
triple of non-negative integers $(n_{\fp,1},n_{\fp,2},n_{\fp,3})$
of distances to the center, as above. We set $N_i = \prod_{\fp\in
A}\fp^{n_{\fp,i}}\in A$ for $i=1,2,3$. Then we have shown that the
modular curve $Y$ is uniquely determined by the triple
$(N_1,N_2,N_3)$. We also see that $p_{i,j}(Y) = Y'_0(N_iN_j)$ for
every pair of coordinates $1\leq i < j \leq 3$.

For $n\geq 4$ such a combinatorial description of pure modular
curves becomes more complicated, as the $\PGL(k_{\fp})$-action is
not $n$-transitive.

\subsection{Modular varieties in $\AA^n$}\label{ModularDefSection}

Let $\pi\in S_n$ be a permutation on $n$ letters, then $\pi$ acts
as a permutation of coordinates on $\AA^n$.

\begin{definition}\label{DrinfeldDefModular1}{\em
An irreducible algebraic variety $X$ in $\AA^n$ is said to be a
{\em modular variety} if it is isomorphic, via some permutation of
coordinates $\pi\in S_n$, to a variety of the form
\begin{equation}\label{Dmodular1}
  \AA^{n_0}\times\prod_{i=1}^{g}Y_i\times\{x\}
\end{equation}
where each $Y_i$ is a pure modular curve in $\AA^{n_i}$ and $x$ is
a CM point in $\AA^{n_{g+1}}$, and $n=n_0+\cdots+n_{g+1}$. The
data
\[
  (\pi, n_0, Y_1,\ldots,Y_g)
\]
is called the {\em type} of $X$.

A reducible variety is modular if all its irreducible components
are modular, and its type is the set of types of the irreducible
components.

A modular variety is {\em pure} if all of the projections $p_i : X
\rightarrow \AA^1$ are dominant on every irreducible component of
$X$. }
\end{definition}

Thus a {\em modular curve} $X$ in $\AA^n$ is a modular variety of
dimension one, i.e. either a pure modular curve, or the product of
a pure modular curve in $\AA^m$ and a CM point in $\AA^{n-m}$.

In the more general case, let $Z=\prod_{i=1}^n X_i$ be a product
of Drinfeld modular curves $X_i\; (= \mbox{compactification of }
\Gamma_i\diagdown\Omega$, where the $\Gamma_i$'s are congruence
subgroups of $\PGLA$). A point $x=(x_1,\ldots,x_n)$ in $Z$ is a CM
point if each $x_i$ corresponds to a CM Drinfeld module (with
$\Gamma_i$-level structure, of course).

A {\em special curve} in $Z$ is the set of points represented by
$(\sigma_1(z),\ldots,\sigma_n(z))\in\Omega^n$ for some
$(\sigma_1,\ldots,\sigma_n)\in\PGLk^n$ and all $z\in\Omega$.

\begin{definition}\label{DrinfeldDefModular2}{\em
An irreducible subvariety $X$ of $Z$ is {\em modular} if there is
a partition $\{1,\ldots,n\} = \coprod_{i=0}^{g+1}S_i$, and $X$ is
given by
\begin{equation}\label{Dmodular2}
 X= \prod_{i\in S_0}X_i\times\prod_{j=1}^{g}Y_j\times\{x\}
\end{equation}
where each $Y_j$ is a special curve in $\prod_{i\in S_j}X_i$ and
$x$ is a CM point in $\prod_{i\in S_{g+1}}X_i$. As before, a
reducible subvariety is modular if all its irreducible components
are modular, and is pure if the projections $p_i : X \rightarrow
X_i$ are dominant for each $i$ from every irreducible component of
$X$. }\end{definition}

But as level structures play no role in these phenomena, any
result concerning $\AA^n$ automatically implies the corresponding
result for $Z=X_1\times\cdots\times X_n$.

\section{Hecke operators}\label{Hecke}
\resetL

Throughout this section, $\fm$ denotes a monic square-free element
of $A$.

\subsection{Hecke operators and Hecke orbits}

\begin{definition}{\em
The {\em Hecke operator} $T_{\AA^n,\fm}$ on $\AA^n$ is the
correspondence given by the image of
\begin{eqnarray*}
Y'_0(\fm)^n & \longrightarrow & \AA^n \times \AA^n \\
\big((x_1,y_1),\ldots,(x_n,y_n)\big) & \longmapsto &
\big((x_1,x_2,\ldots,x_n),(y_1,y_2,\ldots,y_n)\big).
\end{eqnarray*}
We may also view $T_{\AA^n,\fm}$ as a map from subsets of $\AA^n$
to subsets of $\AA^n$, generated by its action on single points:
\begin{eqnarray*}
T_{\AA^n,\fm} : (x_1,\ldots,x_n) \mapsto \{ (y_1,\ldots,y_n)  &| &
\mbox{There exist cyclic isogenies $x_i \rightarrow y_i$} \\
  &  & \mbox{of degree $\fm$ for all $i=1,\ldots,n\}$.}
\end{eqnarray*}
}\end{definition}

We also use the notation $T_{\fm}$ when the $\AA^n$ is clear. We
notice that the operator $T_{\fm}$ is symmetric, in the sense that
$x\in T_{\fm}(y) \Leftrightarrow y\in T_{\fm}(x)$, and that
$T_{\fm}$ is defined over $k$.

Let $X=\cup_{i=1}^r X_i$ be a variety in $\AA^n$, with irreducible
components $X_1,\ldots,X_r$. Then $T_{\fm}(X)$ is also a variety
in $\AA^n$, and $T_{\fm}(X) = \cup_{i=1}^r T_{\fm}(X_i)$.

\begin{definition}{\em
Let $X$ be a variety in $\AA^n$, and suppose all of its
irreducible components have the same dimension. If $X\subset
T_{\AA^n,\fm}(X)$, then we say that $X$ is {\em stabilized} by
$T_{\fm}$, and we define the Hecke operator {\em restricted to}
$X$ by
\[
T_{X,\fm} := \mbox{union of components of
$T_{\AA^n,\fm}\cap(X\times X)$ of maximal dimension.}
\]
}\end{definition}

Whenever we use the notation $T_{X,\fm}$, then it is implicit that
$X$ is stabilized by $T_{\AA^n,\fm}$. The correspondence
$T_{X,\fm}$ is still surjective in the sense that the two
projections $p_X : T_{X,\fm} \rightarrow X$ are surjective. We may
compose Hecke correspondences, and we have the standard property
$T_{\fm_1}\circ T_{\fm_2} = T_{\fm_1\fm_2}$ for $\fm_1,\fm_2\in A$
relatively prime.

\begin{definition}{\em
Let $X\subset\AA^n$ be a variety (possibly $X=\AA^n$), and
$S\subset X$ a subset. Then the {\em Hecke orbit} of $S$ under
$T_{X,\fm}$ is given by
\[
T^{\infty}_{X,\fm}(S) = \cup_{d=1}^{\infty}T^{d}_{X,\fm}(S)
\]
(Here $T^d_{X,\fm}$ means $T_{X,\fm}$ iterated $d$ times.)
}\end{definition}

As there are only finitely many correspondences on the finite set
of irreducible components of $X$, we may decompose $X$ into a
finite disjoint union of Hecke orbits, each orbit being generated
by each of its irreducible components. If $S\subset X_i$ is
Zariski-dense, then $T^{\infty}_{X,\fm}(S)$ is Zariski-dense in
all of $T^{\infty}_{X,\fm}(X_i)$.

Let $x\in\AA^1(\bC)$, and suppose that $x\in T_{\AA^1, \fm}(x)$.
Then $x$ has a cyclic endomorphism, hence is a CM point. For fixed
$\fm \in A$ there are only finitely many such stable points for
$T_{\AA^1,\fm}$, namely the roots of the polynomial
$\Phi_{\fm}(t,t)$. On the other hand, let $x\in\AA^1(\bC)$ be a
given CM point with $\OO=\End(x)$ an order of conductor $f$ in the
CM field $K$. Then for every prime $\fp\in A$, which does not
divide $f$ and decomposes into two principal primes of $K$, we
have $x\in T_{\AA^1, \fp}(x)$. These are precisely the primes
which split completely in the ring class field of $\OO$, hence, by
\v{C}ebotarev, they have density at least $1/2\#\Pic(\OO)$.

\subsection{Some intersection theory}

We define the degree of an irreducible variety $X\subset \AA^n$ of
dimension $d$ as the number of points in the intersection of $X$
with a generic linear subspace of codimension $d$ in $\AA^n$. If
$X$ is not irreducible, then we define its degree to be the sum of
the degrees of its irreducible components of maximal dimension.


We have the following properties, which are easily verified.

\begin{proposition}\label{degrees}
Let $X\subset\AA^n$ be a variety of dimension $d$.
\begin{enumerate}
\item $X$ has at most $\deg(X)$ irreducible components of maximal
dimension. \item {\em (B\'ezout)} If $Y\subset\AA^n$ is another
variety, then $\deg(X\cap Y)\leq\deg(X)\deg(Y)$.
\item $\psi(\fm) \leq \deg(Y'_0(\fm)) \leq 2\psi(\fm)$. \item
$\deg(T_{\AA^n,\fm}(X)) \leq 2^n\psi(\fm)^n\deg(X)$.
\end{enumerate}
\end{proposition}

\begin{proposition}\label{types}
Let $B>0$ and $n\in\NN$ be given. Then there are only finitely
many different types (recall Definition \ref{DrinfeldDefModular1})
of modular varieties $X\subset\AA^n$ with $\deg(X) \leq B$.
\end{proposition}

\paragraph{Proof.}
It suffices to show that there are only finitely many pure modular
curves $Y\subset\AA^n$ with degree less than a given bound. Let
$p_{\{i,i+1\}}(Y)=Y'_0(N_i)$ for $i=1,\ldots, n-1$. Now
$\deg(Y)\geq\deg(p_{\{i,i+1\}}(Y))\geq\psi(N_i)$ for all $i$. But
$\psi(N_i)\rightarrow\infty$ as $N_i$ varies, and the result
follows.\Endproof

\subsection{Preimages in $\Omega^n$}

We have a rigid analytic map $\pi = (j\times\cdots\times j) :
\Omega^n \rightarrow \AA^n(\bC)$. 
For each irreducible component
$X_i$ of $X$ we choose an irreducible component $Z_i$ of the rigid
analytic variety $\pi^{-1}(X_i)\subset\Omega^n$. We set $Z=\cup_i
Z_i$. The group $\PGLki^n$ acts on $\Omega^n$ and the
$\PGLA^n$-orbit of $Z_i$ is all of $\pi^{-1}(X_i)$. We want to
describe the Hecke operators acting in the space $\Omega^n$.

For the following discussion of matrices, see \cite{Bae92}. We let
$\Delta^*_{\fm}$ denote the set of $2\times 2$ matrices over $A$
with determinant in $\Fq^*\fm$ and whose entries have no factor in
common. Then $\GLA$ acts from the right on $\Delta^*_{\fm}$, and
the representatives of $\Delta^*_{\fm}/\GLA$ may be chosen of the
form
\[
t_i = \matrix{a_i}{b_i}{0}{d_i}, \qquad i\in
\II:=\{1,\ldots,\psi(\fm)\}
\]
where $a_i,d_i$ are monic, $a_id_i = \fm$ and $|b_i| < |a_i|$. For
$i=(i_1,\ldots,i_n)$ ranging through $\II^n$ we denote by $t_i =
(t_{i_1},\ldots,t_{i_n})$ the resulting representatives of
$(\Delta^*_{\fm})^n/\GLA^n$.

For each $Z_i$ we define $\JJ_{Z_i} \subset \II^n$ as the set of those 
indices $j$ for which $t_j(Z_i) \subset \pi^{-1}(X)$. 
Let $x\in X_i(\bC)$, and choose some $z\in Z_i$ with $\pi(z)=x$.
Then the action of $T_{\AA^n,\fm}$ and $T_{X,\fm}$ on $x$ are given
by $T_{\AA^n,\fm}(x) = \{\pi(t_j(z)) \;|\; j\in\II^n\}$ and 
$T_{X,\fm}(x) = \{\pi(t_j(z)) \;|\; j\in\JJ_{Z_i}\}$. In particular,
$\JJ_{Z_i}$ is non-empty, as $T_{X,\fm}$ is a surjective correspondence.

%

\subsection{Surjectivity of projections}\label{surjsection}

%
%

\begin{theorem}\label{surjectivity}
Let $X\subset \AA^n$ be a variety all of whose irreducible
components have the same dimension, and suppose that $X \subset
T_{\fm}(X)$ for some square-free $\fm\in A$ which is a product of
distinct primes $\fp\in A$ of even degree satisfying $|\fp|\geq
\max(13,\deg X)$. 
Let $X_i$ be an irreducible component of $X$ for which the projection
$p_I : X_i \rightarrow \AA^I$ is dominant.
Then the projection
\begin{equation}\label{eqsurj1}
p_I : \JJ_{Z_i} \longrightarrow \II^I
\end{equation}
is surjective. In particular, let $x\in X_i(\bC)$. Then the
projection of finite sets
\begin{equation}\label{eqsurj2}
p_I : T_{X,\fm}(x) \longrightarrow T_{\AA^I,\fm}(p_I(x))
\end{equation}
is surjective.
\end{theorem}

\paragraph{Proof.}

Clearly, the surjectivity of (\ref{eqsurj2}) follows from the
surjectivity of (\ref{eqsurj1}), which in turn follows from the
surjectivity of (\ref{eqsurj2}) for a {\em generic} point $x\in
X_i$.

So we suppose $x\in X_i$ is generic. Denote by $T_{X,\fm,i} =
T_{X,\fm}\cap(X_i \times X)$ the restriction of the source of the
Hecke correspondence $T_{X,\fm}$ to the component $X_i$. Consider
the following diagram
\[\xymatrix{
T_{X,\fm,i}\ar[rr]^{p_I\times p_I}\ar[dd]^{p_{X_i}}\ar[dr]^{f_i} 
& & T_{\AA^I,\fm}\ar[dd]^{p_{\AA^I}} \\
& X_i \times_{\AA^I} T_{\AA^I,\fm}\ar[ur]\ar[dl] & \\
X_i\ar[rr]^{p_I} & & \AA^I, } \] where the vertical arrows are
projections onto the sources of the respective correspondences,
and the horizontal arrows are projections onto the coordinates in
$I$. There exists a canonical map $f_i$ from $T_{X,\fm,i}$ to the
fibered product $X_i \times_{\AA^I} T_{\AA^I,\fm}$, which is
generically finite, as $p_{X_i} : T_{X,\fm,i} \rightarrow X_i$ is
generically finite. Clearly $X_i \times_{\AA^I} T_{\AA^I,\fm}$ has
the same dimension as $T_{X,\fm}$ and $X$, and it is irreducible,
as we will show below. It follows that $f_i$ is dominant.

Now let $x_I = p_I(x)$, and let $y_I \in T_{\AA^I,\fm}(x_I)$. Then
$\big(x,(x_I,y_I)\big)$ is a generic point on the fibered product,
hence has a preimage $(x,y)$ under $f_i$. We see that $y\in
T_{X,\fm,i}(x)=T_{X,\fm}(x)$ and $p_I(y)=y_I$, so we have shown
that (\ref{eqsurj2}) is surjective.

It remains to show that $X_i \times_{\AA^I} T_{\AA^I,\fm}$ is
irreducible. This will follow if the function fields of $X_i$ and
$T_{\AA^I,\fm}\cong (Y'_0(\fm))^I$ over $\bC$ are linearly
disjoint over the function field of $\AA^I$ over $\bC$. Recall
from Proposition \ref{GaloisCover} that the modular curve
$Y_2(\fm)$ covers $Y_0(\fm)$ and is Galois over $Y(1)=\AA^1$ with
Galois group
\[
\Gal(Y_2(\fm)/Y(1)) \cong \PSL(A/\fm A) \cong
\prod_{\fp|\fm}\PSL(A/\fp A).
\]
On the other hand, let $L$ be an intermediate field $\bC(\AA^I)
\subset L \subset \bC(X_i)$ which is purely transcendental over
$\bC(\AA^I)$ and for which $[\bC(X_i):L] \leq\deg(X_i)\leq\deg(X)$
is finite.
Then $L\cap\bC(Y_2(\fm)^I) = \bC(\AA^I)$, and $\bC(Y_2(\fm)^I)$ is
Galois over $\bC(\AA^I)$,
so it follows
that $L$ and $\bC(Y_2(\fm)^I)$ are linearly disjoint over
$\bC(\AA^I)$. Denote by $L_{\fm}$ the field $L
\otimes_{\bC(\AA^I)} \bC(Y_2(\fm)^I)$. Then we have
\begin{eqnarray*}
\bC(X_i)\otimes_{\bC(\AA^I)}\bC(T_{\AA^I,\fm})
& \subset & \bC(X_i)\otimes_{\bC(\AA^I)}\bC(Y_2(\fm)^I) \\
& = & \bC(X_i)\otimes_L L_{\fm}.
\end{eqnarray*}
But now $L_{\fm}$ is Galois over $L$, with group
$\Gal(L_{\fm}/L)\cong\PSL(A/\fm A)^I$. Moreover, as $|\fp|\geq
13$ for all $\fp|\fm$, it follows that this group has no proper
subgroup of index less than $|\fp|+1$ (see for example
\cite{Huppert}, Hauptsatz 8.27 and Satz 8.28).

On the other hand, $[\bC(X_i):L] \leq\deg(X) < |\fp|+1$, so
$\bC(X_i)\cap L_{\fm}=\bC(\AA^I)$. It follows that $\bC(X_i)$ and
$L_{\fm}$ are linearly disjoint, hence
$X_i \times_{\AA^I} T_{\AA^I,\fm}$ is irreducible, as required.
\Endproof

For the next two corollaries, we assume $X\subset\AA^n$ is a
variety, with irreducible components $X_i,\;i=1,\ldots,r$, which
are all of the same dimension. We assume further that $X \subset
T_{\AA^n,\fm}(X)$ for some square-free $\fm\in A$, composed of
distinct primes $\fp\in A$, each of even degree and satisfying
$|\fp| \geq \max(13,\deg X)$.

\begin{corollary}\label{SurjCor1}
Suppose that the projection $p_1 : X_i \rightarrow \AA^1$ onto the
first coordinate is dominant for all $i=1,\ldots,r$. Let
$x_1\in\AA^1$ such that $x_1 \in T_{\AA^1,\fm}(x_1)$. Let $X_{x_1}
= X \cap (\{x_1\}\times\AA^{n-1})$. Then
\[
X_{x_1} \subset T_{X,\fm}(X_{x_1}).
\]
\end{corollary}

\paragraph{Proof.}
Let $x\in X_{x_1}$. Then setting $I=\{1\}$ in Theorem
\ref{surjectivity}, we see that
\[
p_1 : T_{X,\fm}(x) \longrightarrow T_{\AA^1,\fm}(x_1)
\]
is surjective. Let $y\in T_{X,\fm}(x)$ be a preimage of $x_1\in
T_{\AA^1,\fm}(x_1)$. Then $y\in X_{x_1}$ and $x\in T_{X,\fm}(y)$,
hence $x\in T_{X,\fm}(X_{x_1})$, as required. \Endproof

\begin{corollary}\label{SurjCor2}
Let $x\in X_i$. Then the Hecke orbit $T^{\infty}_{X,\fm}(x)$ is
Zariski-dense in the Hecke orbit $T^{\infty}_{X,\fm}(X_i)$.
\end{corollary}

\paragraph{Proof.}
Clearly, we may suppose that $\dim(X)\geq 1$.
Let $I\subset\{1,\ldots,n\}$ be such that $\#I=\dim(X)$ and the
projection $p_I : X_i \rightarrow \AA^I$ is dominant. 

We claim that $p_I : X_j \rightarrow \AA^I$ is also dominant for 
every irreducible component $X_j$ of $T^{\infty}_{X,\fm}(X_i)$. By induction, 
it suffices to prove the claim for $X_j \subset T_{X,\fm}(X_i)$.

Let $x_I \in \AA^I$ be a generic point. Then there is some $x\in X_i$ with 
$p_I(x)=x_I$. At least one point $y\in T_{X,\fm}(x)$ lies on $X_j$, and 
$p_I(y)=y_I \in T_{\AA^I,\fm}(x_I)$.
So it follows that every generic $x_I\in \AA^I$ is $\fm$-isogenous to some $y_I$ 
coming from $X_j$, in other words, $T_{\AA^I,\fm}(p_I(X_j))$ is Zariski-dense 
in $\AA^I$. It follows that 
$\dim(p_I(X_j))=\dim(T_{\AA^I,\fm}(p_I(X_j))) = \dim(\AA^I)$, 
and so $p_I(X_j)$ is Zariski-dense in $\AA^I$, which proves the claim.

Now we apply Theorem
\ref{surjectivity} to obtain a surjection
\[
p_I : T^{\infty}_{X,\fm}(x) \longrightarrow
T^{\infty}_{\AA^I,\fm}(p_I(x)).
\]
This last set is Zariski-dense in $\AA^I$, as
$T^{\infty}_{\AA^I,\fm}(p_I(x)) = \prod_{j\in
I}T^{\infty}_{\AA^1,\fm}(x_j)$ is a product of infinite subsets of
$\AA^1$. As the projection $p_I : X \rightarrow \AA^I$ is
generically finite, it follows that $T^{\infty}_{X,\fm}(x)$ must
be Zariski-dense on at least one component $X_j$ of
$T^{\infty}_{X,\fm}(X_i)$, hence on all of
$T^{\infty}_{X,\fm}(X_i)$. \Endproof

\begin{remark}{\em
As $T_{\fm}$ is defined over $k$, we may replace the word
``irreducible'' by ``$F$-irreducible'' everywhere in the preceding
sections, for any field $F \supset k$ over which the relevant
varieties are defined. In particular, it follows from Corollary
\ref{SurjCor2} above, that if $X$ is a variety defined over $F$,
$X_i$ is an $F$-irreducible component of $X$, and $x\in X_i$, then
the Hecke orbit $T^{\infty}_{X,\fm}(x)$ is Zariski-dense on $X_i$.
}\end{remark}

\subsection{Curves stabilized by Hecke operators}\label{TopologySection}

We are now ready to prove a fundamental result: a characterization
of the modular curves $Y'_0(N)$ in terms of Hecke operators.

\begin{theorem}\label{Topology}
Let $X\subset\AA^2$ be an irreducible algebraic curve, and suppose
$X\subset T_{\AA^2,\fm}(X)$ for some square-free $\fm\in A$,
$|\fm|>1$, composed of primes $\fp\in A$ of even degree satisfying
$|\fp|\geq\max(13,\deg X)$. Then $X=Y'_0(N)$ for some $N\in A$.
\end{theorem}

The proof will occupy the next three sections.

If $X=\{x\}\times\AA^1$ or $X=\AA^1\times\{x\}$, then $x$ is a CM
point (as it is stabilized by $T_{\AA^1,\fm}$), and so $X$ is
modular. So we may assume that the projections $p_i : X
\rightarrow \AA^1$ are dominant, and have degree $1 \leq d_i \leq
\deg(X)$, for $i=1,2$.

The group $G:=\PGLki^2$ acts on $\Omega^2$, and we also define the
following groups: $S:=\PSLki^2,\; \Gamma:=\PGLA^2$, and
$\Sigma:=\PSLA^2$. As before, we choose an irreducible component
$Z$ of the rigid analytic variety $\pi^{-1}(X)$. Let $G_Z$ be the
stabilizer of $Z$ under the action of $G$, it is a closed analytic
subgroup of $G$. We also define $S_Z:=G_Z\cap S,\; \Gamma_Z :=
G_Z\cap\Gamma$, and $\Sigma_Z := G_Z\cap \Sigma$. Our aim is to
investigate the structure of $S_Z$, under the hypothesis that
$X\subset T_{\fm}(X)$, and hence conclude that $X$ must be a
modular curve.

So our whole approach is similar to that of Edixhoven
\cite{Edixhoven98}, but with slightly different details, for
example the action of $G$ on $\Omega^2$ is not transitive, the
topology is ultrametric, and Lie theory works a bit differently in
characteristic $p$, so we replace it by explicit calculations.

We denote by $pr_i : G \rightarrow \PGLki$ the two projections,
$i=1,2$. The following lemma holds for an arbitrary curve $X$
(with non-constant projections).

\begin{lemma}\label{projlemma}
$\;$
\begin{enumerate}
\item The two projections $pr_i:G_Z \rightarrow \PGLki$ are
injective. \item $pr_i(\Gamma_Z)$ has index at most $d_i$ in
$\PGLA$.
\end{enumerate}
\end{lemma}

\paragraph{Proof. (1)}
Let $K = \ker(pr_2 : G_Z \rightarrow \PGLki)$. Then $K$ is the
stabilizer of $Z$ in $\PGLki\times\{1\}$, and stabilizes $Z_z =
Z\cap(\Omega\times\{z\})$, for any $z\in\Omega$. But $Z_z$ is
discrete, and we may choose $z$ in such a way that $Z_z$ contains
a non-quadratic element, whose stabilizer is trivial, so it
follows that $K$ is discrete. Now $K \triangleleft
\PGLki\times\{1\}$, which has no non-trivial discrete normal
subgroups, thus $K=\{1\}$. The same holds for the other
projection.

\subparagraph{(2)} We factor the map $\pi$ as follows:
\[\xymatrix@R=2pt{
\Omega\times\Omega \ar[r]^{\pi_1} & \AA^1\times\Omega
\ar[r]^{\pi_2}
  & \AA^1\times\AA^1 \\
Z \ar[r] & W \ar[r] & X }\]

Here $W=\pi_1(Z)$ is an irreducible component of
$Y=\pi_2^{-1}(X)=\PGLA\cdot\nolinebreak W$. Let $S$ be the set of
all $x$'s in $X$ for which every (equivalently at least one) point
of $\pi^{-1}(x)$ lies in more than one component of $\pi^{-1}(X)$.
Then $S$ lies in the finite set consisting of the singular points
of $X$ as well as those with at least one coordinate equal to $0$.

Let $X'=X-S$, and let $Z'$ and $W'$ be the corresponding
preimages, then the map $\pi : Z' \rightarrow X'$ is a quotient
for the action of $\Gamma_Z$, and $\pi_2 : W' \rightarrow X'$ is a
quotient for the action of $pr_2(\Gamma_Z)$. It follows that
$pr_2(\Gamma_Z)$ is the stabilizer of $W$ for the action of
$\PGLA$, hence the irreducible components of $Y$ correspond to the
cosets $\PGLA / pr_2(\Gamma_Z)$, so the index is the number of
these components.

On the other hand, $Y=\pi_2^{-1}(X)$ is the fibered product of the
maps $p_2 : X \rightarrow \AA^1$ and $j : \Omega \rightarrow
\AA^1$, hence it has at most $d_2$ irreducible components. Again,
the same holds for the other projection. \Endproof

\subsection{The structure of $S_Z$}\label{matrixsection}

Now we make use of the fact that $X\subset T_{\fm}(X)$. Let
$i\in\JJ_Z$, then $t_i(Z)\subset\pi^{-1}(X)$ (by definition;  
remember that $X$ is irreducible), and
thus there is some $\gamma_i\in\Gamma$ such that $g_i:=\gamma_i
t_i \in G_Z$. Moreover, applying Theorem \ref{surjectivity} to the
projection $p_i : X \rightarrow \AA^1$ we see that the projection
$p_i : \JJ_Z \rightarrow \II=\{1,\ldots,\psi(\fm)\}$ is
surjective. This gives us many non-trivial elements in $G_Z$. More
precisely, we will study the projections $H_i = pr_i(G_Z)$, and
show that they each contain $\PSLki$.

From Lemma \ref{projlemma} follows that $\PGLA\cap H_1$ has finite
index in $\PGLA$. Let $R$ be a finite set of representatives of
$\PGLA/(\PGLA\cap H_1)$. The group $\GL(A)$ acts from the right on
the set of left cosets $\GL(A)\diagdown\Delta^*_{\fm}$. We claim
that for any string $i_1\ldots i_n$ of elements in $\II$, and any
$a \in \GL(A)$, we can construct an element of the form $\gamma
t_{i_n}t_{i_{n-1}}\cdots t_{i_1}a$ in $H_1$, for some $\gamma\in
R$ depending on the string and on $a$. Indeed, by induction we
need only show that, given $a_1 \in \GL(A)$ and $i_1 \in \II$, we
can construct an element of the form $\gamma_1 t_{i_1} a_1$ in
$H_1$. This element is constructed as follows. Let $a_1$ act from
the right on the coset $\GL(A)\cdot t_{i_1}$, to obtain another
coset $\GL(A)\cdot t_{i_1} a_1 = \GL(A)\cdot t_j$. Then $t_{i_1}
a_1 = \gamma'_j t_j$, and multiplying on the left with a suitable
element $\gamma_1$ of $R$ gives $\gamma_1 t_{i_1} a_1 =
\gamma'_1\gamma_jt_j = \gamma'_1g_j\in H_1$, with $\gamma'_1\in
H_1\cap\PGLA$. This proves the claim.

Multiplying by a suitable power of the scalar $\fm$, we see that
for any $x\in A[1/\fm]$ and any $a\in\GL(A)$, there exists
$\gamma_{x,a}\in R$ such that $\gamma_{x,a}\matrix{1}{x}{0}{1} a
\in H_1$.

The group $\PSL(A[1/\fm])$ is generated by $\PSLA$ and elements of
the form $\matrix{1}{x}{0}{1}$, hence for any $g\in
\PSL(A[1/\fm])$, we can construct an element $\gamma_g g \in H_1$,
for some $\gamma_g \in R$, obtained by multiplying together
suitable elements of the form $\gamma_{x,a}\matrix{1}{x}{0}{1} a
\in H_1$. It follows that $H_1\cap\PSL(A[1/\fm])$ has finite index
in $\PSL(A[1/\fm])$, which is dense in $\PSLki$.

\begin{lemma}
$G_Z$ is not discrete
\end{lemma}

\paragraph{Proof.}
Assume that $G_Z$ is discrete. Choose a non-quadratic point
$z=(z_1,z_2)$ in $Z$. Then its orbit $G_Z\cdot z$ is discrete in
$Z$, so $\pi(G_Z\cdot z)$ is discrete in $X$, as $\Gamma_Z \subset
G_Z$. Next, $p_1\big(\pi(G_Z\cdot z)\big)$ is discrete in $\AA^1$
(as $p_1 : X \rightarrow \AA^1$ is finite), and thus
$j^{-1}\Big(p_1\big(\pi(G_Z\cdot z)\big)\Big)$ is discrete in
$\Omega$. But from above we see that this set contains the orbit
$\Big(H_1\cap\PSL(A[1/\fm])\Big)\cdot z_1$, which is not discrete.
This is a contradiction. \Endproof

So we see that $G_Z$ is a closed analytic subgroup of $G$ which is
not discrete, hence of dimension at least one. The projection
$pr_1 : S_Z \rightarrow \PGLki$ is injective, so we see that $H_1
= A_1 \smallsetminus B_1$, where $A_1, B_1$ are analytic sets and
$\dim(A_1)=\dim(G_Z) > \dim(B_1)$. Thus, there exists a point
$x\in H_1$ and a closed (in the non-archimedean topology)
neighborhood $A_x$ of $x$ in $\PGLki$ such that $A_x\cap H_1$ is
closed in $\PGLki$, i.e. $H_1$ is locally closed at $x$. Since
$H_1$ is a topological group, it is closed in $\PGLki$.


Now, $\PSL(A[1/\fm])$ is dense in $\PSLki$, so $H_1\cap\PSLki$ has
finite index in $\PSLki$, which is simple, hence $\PSLki \subset
H_1$. In particular, $H_1$ has finite index in $\PGLki$. Of
course, the same holds for $H_2 = pr_2(G_Z)$.

Goursat's lemma says that $G_Z$ is of the form
\[
G_Z = \{ \big(g,\rho(g)\big) \;|\; g \in H_1\}
\]
for some isomorphism $\rho : H_1 \stackrel{\sim}{\rightarrow} H_2
\subset \PGLki$. Now, as $\PSLki$ is simple and the image of
$\rho$ has finite index in $\PGLki$, it follows that $\rho$
restricts to an automorphism on $\PSLki$. Thus we have shown
\[
G_Z \cap(\PSLki)^2 = S_Z = \{ (g,\rho(g)) \;|\; g \in \PSLki\}
\]
for some $\rho\in\Aut(\PSLki)$.

Every automorphism of $\PSLki$ is of the form $g \mapsto
hg^{\sigma}h^{-1}$ for some $h\in \PGLki$ and
$\sigma\in\Aut(\kinf)$, see \cite{Hua}.

By the definition of $\Sigma_Z$ and the structure of $S_Z$, we see
that $h\cdot pr_1(\Sigma_Z)^{\sigma} \cdot h^{-1} \subset \PSLA$.
On the other hand, Lemma \ref{projlemma} tells us that
$pr_1(\Sigma_Z)$ has finite index in $\PSLA$. This in turn
severely restricts $h$ and $\sigma$:

\begin{proposition}\label{prop1}
Let $G$ be a subgroup of finite index in $\PSL(A)$, and suppose
that \linebreak $hG^{\sigma}h^{-1} \subset \PGLk$, for some $h\in\PGLki$ and
$\sigma\in\Aut(\kinf)$. Then $h\in\PGLk$ and $\sigma(T) = uT+v$
for some $u\in\Fq^*,\; v\in\Fq$, and $\sigma(\Fq)=\Fq$.
\end{proposition}

\subparagraph{Proof.} Firstly, let $h=\matrix{a}{b}{c}{d}$ and let
$r=\det(h)$. As $\kinf^*/\kinf^{*2}$ may be represented by
$\{1,\alpha,T,\alpha T\}$, for some non-square $\alpha\in\Fq$, we
may assume that $r\in k$.

Denote by $B_1 = \matrix{1}{*}{0}{1}$ and $B_2 =
\matrix{1}{0}{\mbox{$*$}}{1}$ the two Borel subgroups of $\PSLA$.
The group $G$ has finite index in $\PSLA$, so it
follows that $G\cap B_1$ and $G\cap B_2$ are of finite index in
$B_1$ and $B_2$, respectively. Hence
\[
A_0^+ := \{x \in A \;|\; \matrix{1}{x}{0}{1}, \matrix{1}{0}{x}{1}
\in G\}
\]
has finite index in the additive group $A^+$ of $A$.
%
Now for every $x\in A_0^+$ we have
\[
\matrix{a}{b}{c}{d}\matrix{1}{x}{0}{1}^{\sigma}\matrix{a}{b}{c}{d}^{-1}
= \matrix{1 - \frac{ac}{r}\sigma(x)}{\frac{a^2}{r}\sigma(x)}
{-\frac{c^2}{r}\sigma(x)}{1+\frac{ac}{r}\sigma(x)} \in \PGLk,
\]
and it follows that
\begin{equation}\label{eq1}
ac\sigma(x),\; a^2\sigma(x),\; c^2\sigma(x) \in k.
\end{equation}
Likewise, from
$\matrix{a}{b}{c}{d}\matrix{1}{0}{x}{1}^{\sigma}\matrix{a}{b}{c}{d}^{-1}
\in \PGLk$ follows
\begin{equation}\label{eq2}
bd\sigma(x),\; b^2\sigma(x),\; d^2\sigma(x) \in k.
\end{equation}
From (\ref{eq1}), (\ref{eq2}) and $ad-bc\in k$, we may deduce, in that order,
\begin{equation}\label{eqlist}
\frac{a}{c}, \frac{b}{d}, cd, ab, \frac{a}{b}\sigma(x),
\frac{c}{d}\sigma(x), \frac{c^2}{d^2} \in k
\end{equation}
from which follows that
\begin{equation}\label{eq6}
\sigma(x)^2 \in k \qquad \forall x\in A_0^+.
\end{equation}
Now as $A_0^+$ has finite index in $A^+$ it follows that there
exists a pair $x_1\neq x_2 \in A_0^+$ such that $y = x_1^2 - x_2^2
\in A_0^+$. We have $\sigma(y) = \sigma(x_1)^2 - \sigma(x_2)^2 \in
k$. Substituting this $y$ for $x$ in (\ref{eqlist}) shows that
in fact
\begin{equation}\label{eq7}
\frac{a}{b} \in k,\; \frac{c}{d} \in k \quad\mbox{and}\quad
\frac{a}{d}=\frac{a}{b}\frac{b}{d} \in k.
\end{equation}
It follows firstly that $h\in\PGLk$, and secondly that $\sigma(x)
\in k$ for all $x\in A_0^+$. As the elements of $A_0^+$ generate
$k$ as a ring, $\sigma(x) \in k$ for all $x\in k$. It remains to
characterize those automorphisms $\sigma$ for which
$\sigma(k)\subset k$.

Let $R=\Fq[[1/T]]=\{x \in \kinf \;|\; |x| \leq 1\}$. Then $R$ is
the unique valuation ring of $\kinf$. It is characterized by the
property: $x\in R$ or $x^{-1}\in R$ for all $x\in\kinf$ and $R\neq
\kinf$. This property must be preserved by $\sigma$, so $\sigma(R)
\subset R$. So $\sigma$ also preserves $k\cap R = A$, and the only
automorphisms that send polynomials to polynomials are of the form
$\sigma(T) = uT + v$, for some $u\in\Fq^*,\; v\in\Fq$, and
$\sigma(\Fq) = \Fq$. \Endproof

\subsection{Completing the proof of Theorem \ref{Topology}}

\paragraph{Proof of Theorem \ref{Topology}.}

From Proposition \ref{prop1} follows that
\[
S_Z = \{ (g,hg^{\sigma}h^{-1}) \;|\; g \in \PSLki\},
\]
where $h\in\PGLk,\;\sigma(T) = uT+v$ and $\sigma(\Fq)=\Fq$. There
is some $t\in\NN$ such that $\sigma(\alpha) = \alpha^{p^t}$ for
all $\alpha\in\Fq$, as $\sigma|_{\Fq} \in \Gal(\Fq/\Fp)$.

We let $f=(T^q - T)^{q-1}$, then $\sigma(f) = f$. Let $F =
\Fp((1/f))$. This is a complete subfield of $\kinf$ and $\sigma$
acts trivially on $F$.

Now fix some non-square $\alpha\in\Fq$, and define the set
\[
P = \{ z \in \Omega \;|\; z^2 = \alpha e,\; e \in F \}.
\]
This is an uncountable subset of $\Omega = \bC \smallsetminus
\kinf$, as $\sqrt{\alpha}\not\in\kinf$.

Next, we notice that $\sigma(\alpha e) = \alpha^{p^t} e = \beta^2
\alpha e$, where we set $\beta = \alpha^{(p^t - 1)/2} \in \Fq^*$
(remember that $p$ is odd).

Let $z_1 = \sqrt{\alpha e} \in P$ and
\[
S_1 = \Stab_{\PSL(F)}(z_1) = \left\{\matrix{a}{b}{c}{d} \;|\;
a=d,\; b=c\alpha e,\; ad-bc = 1\right\} \diagup \{\pm 1\},
\]
which is a one-dimensional Lie-group over $F$.

Now let $z_2 \in \Omega$ such that $(z_1,z_2) \in Z$, and consider
the ``$S_1$-orbit'' of $(z_1,z_2)$:
\[
\{\big(g(z_1),hg^{\sigma}h^{-1}(z_2)\big) \;|\; g \in S_1\}
\subset Z\cap(\{z_1\}\times\Omega).
\]
This set is discrete, but the group $S_1$ is not, hence there
exists some non-trivial $g\in S_1$ such that $g$ fixes $z_1$ (by
definition of $S_1$) and $hg^{\sigma}h^{-1}$ fixes $z_2$. But
$g^{\sigma}$ fixes the point $z_1^{\sigma} := \sqrt{\sigma(\alpha
e)} = \beta z_1$, so we see that $hg^{\sigma}h^{-1}$ fixes both
$z_2$ and $h(\beta z_1) = h'(z_1)$, where we have written $h' =
h\circ\matrix{\beta}{0}{0}{1}\in\PGLk$. However, any non-trivial
element of $\PGLki$ fixes at most two points of $\Omega$, namely a
conjugate pair of quadratic points. So $z_2$ and $h'(z_1)$ are
conjugate. So we get either $z_2 = h'(z_1)$ or $z_2 = h'(-z_1)$.
As $j(z_1) = j(-z_1)$, we get either $\big(j(z_1),j(h'(z_1))\big)$
or $\big(j(-z_1),j(h'(-z_1))\big)$ on the curve $X$ in
$\AA^2(\bC)$.

Let $N$ be the degree of $h'$. Then we see that the points
$\big(j(z_1),j(h'(z_1))\big)$ and \linebreak $\big(j(-z_1),j(h'(-z_1))\big)$
lie on $Y'_0(N)$ (which is independent of $z_1$). We get such a
point for each $z_1 \in P$, and $P$ is uncountable whereas the
fibers of $j$ are countable, so it follows that $X(\bC)\cap
Y'_0(N)(\bC)$ is infinite, hence $X=Y'_0(N)$. This completes the
proof of Theorem \ref{Topology}. \Endproof

By considering various projections onto pairs of coordinates, we
immediately get

\begin{corollary}\label{ToTopology}
Let $X\subset\AA^n$ be an irreducible algebraic curve, and suppose
that $X\subset T_{\AA^n,\fm}(X)$ for some square-free $\fm\in A$,
$|\fm|>1$, and composed of primes $\fp$ of even degree and
satisfying $|\fp| \geq \max(13,\deg X)$. Then $X$ is a modular
curve.
\end{corollary}

\subsection{Varieties stabilized by Hecke operators}\label{TopologySection2}

In this section we generalize Theorem \ref{Topology} to
subvarieties of higher dimensions.

\begin{theorem}\label{Topology2}
Let $F$ be a field lying between $k$ and $\bC$. Let $X\subset
\AA^n$ be an $F$-irreducible variety, containing a CM point $x\in
X(\bC)$. Suppose that $X\subset T_{\AA^n,\fm}(X)$ where $\fm\in
A$, $|\fm|>1$ is square-free, composed of primes $\fp$ of even
degree and satisfying $|\fp|\geq \max(13,\deg X)$. Then $X$ is a
modular variety.
\end{theorem}

\paragraph{Proof.} We know from Corollary \ref{SurjCor2}, and the subsequent Remark,
that the Hecke orbit $S=T^{\infty}_{X,\fm}(x)$ is Zariski-dense in
$X$. In particular, it is Zariski-dense on every (geometrically)
irreducible component. So now we assume that $X$ is geometrically
irreducible (but not necessarily stabilized by $T_{\AA^n,\fm}$ -
indeed, all we need is a dense Hecke orbit of CM points). All the
points in $S$ are CM points, isogenous coordinate-wise to $x$. As
CM points are defined over $k^{sep}$, so is $X$. So we may also
assume that $X$ is defined over a finite Galois extension (again
denoted $F$) of $k$.

\subparagraph{Step 1.} Write $x=(x_1,\ldots,x_n)$, and let
$\OO_i=\End(x_i)$ be an order of conductor $f_i$ in the imaginary
quadratic field $K_i$, for each $i=1,\ldots,n$. Set $K=K_1\cdots
K_n$ and $f=f_1\cdots f_n$, and define
\[
\cP = \{l\in A \;|\; \mbox{monic prime, of even degree, split
completely in $FK$ and $l\nmid f\fm$}\}.
\]
This set is infinite (\v{C}ebotarev).

Let $x'=(x'_1,\ldots,x'_n)\in S$, then each $\OO'_i=\End(x'_i)$ is
an order of conductor $f'_i$ in $K_i$ (the CM fields are the same,
as $x_i$ and $x'_i$ are isogenous). Furthermore, all the prime
factors of $f'_i$ are factors of $f_i$ and of $\fm$. It follows
that every $l\in\cP$ splits also in $\OO'_i$. Set
$M=K(x'_1,\ldots,x'_n)$ and let $\fL$ be a prime of $FM$ lying
over $l$. Denote by $\fL_i$ the restriction of $\fL$ to the field
$K_i(x'_i)$. From CM Theory (Theorem \ref{CMtheory}) follows that
$l$ is unramified in $M$, hence also in $FM$. Let $\sigma =
(\fL,FM/k)$ be the Frobenius element. Set $\sigma_i =
\sigma|_{K_i(x'_i)} = (\fL_i,K_i(x'_i)/k)$. As $l$ splits in
$K_i$, we have in fact $\sigma_i = (\fL_i,K_i(x'_i)/K_i)$. CM
theory then tells us that there is a cyclic isogeny $x'_i
\rightarrow \sigma_i(x'_i)$ of degree $l$. Now $\sigma$ fixes $F$,
and we have
\[
x' \in X\cap T_{\AA^n,l}(X^{\sigma}) = X\cap T_{\AA^n,l}(X).
\]
This holds for every $x'$ in the Zariski-dense set $S$, so it
follows that
\begin{equation} \label{HeckeEq1}
X \subset T_{\AA^n,l}(X).
\end{equation}
Moreover, (\ref{HeckeEq1}) holds for every $l\in\cP$.

\subparagraph{Step 2.} Now we use induction on $d=\dim(X)$, and
suppose $d\geq 2$.

We may assume without loss of generality that the projection $p_1
: X \rightarrow \AA^1$ is dominant. Now we may choose an infinite
subset $\{x^1,x^2,\ldots,\}\subset S$ of points, written $x^j =
(x^j_1,\ldots,x^j_n)$, such that the first coordinates $x^j_1$ are
distinct, for $j\in\NN$. For each $j$ we may find $l_j \in \cP$
such that $x^j_1 \in T_{\AA^1,l_j}(x^j_1)$ and $|l_j|\geq
\max(13,\deg X)$. In fact, $\cP$ contains infinitely many such
primes, namely those which split completely in the ring class
field of $\End(x^j_1)$.

For each $j\in\NN$ we consider the ``slice''
\[
X_j = X\cap(\{x^j_1\}\times\AA^{n-1}),
\]
which satisfies $X_j \subset T_{X,l_j}(X_j)$ (Corollary
\ref{SurjCor1}), $\dim(X_j)=d-1$ and $x^j \in X_j$. Let $X'_j$ be
an irreducible component of $X_j$ containing $x^j$. Then the Hecke
orbit $T^{\infty}_{X_j,l_j}(x^j)$ is Zariski-dense in $X'_j$. As
in Step 1 above, we can find infinitely many primes $\fp$ such
that $T_{\AA^n,\fp}$ stabilizes $X'_j$, so from the induction
hypothesis follows that $X'_j$ is modular.

Now $\deg(X'_j) \leq \deg(X)$, and there are only finitely many
types of modular varieties of bounded degree (Proposition
\ref{types}), so it follows that we have an infinite subset
$I\subset \NN$ and some $\pi \in S_n$ such that, after permutation
of coordinates by $\pi$,
\[
X'_j = Y \times \{y_j\} \quad\forall j\in I,
\]
where $Y\subset \AA^{n-m}$ is a fixed modular variety, and $y_j
\in \AA^m$ is a CM point, for some $m\geq 1$. 
Let $Y'\subset\AA^m$ be the Zariski-closure of $\{y_j \;|\; j\in I\}$,
then $\dim(Y')\geq 1$. Now the Zariski-closure of 
$\{X'_j \;|\; j\in I\}$ is equal to $Y\times Y'$, is contained in
$X$ and has dimension at least $\dim(Y)+1 = \dim(X)$. It follows 
that $X=Y\times Y'$, with $Y'$ an irreducible curve.
%
%

Moreover, $Y'$ is stabilized by
the Hecke operators $T_{\AA^m,l}$ for all $l\in\cP$, hence is
itself modular. It follows that $X$ is modular, which is what we
set out to prove. \Endproof

\section{Heights of CM points}\label{Heights}
\resetL

\subsection{Estimating class numbers}\label{sectionClassNumbers}

We now want to derive a lower bound for the class number of an
order in an imaginary quadratic function field. Our standard
reference to facts about function fields is \cite{Stichtenoth93}.

Let $F$ be a global function field of genus $g$ and exact field of
constants $\Fq$, and denote by $h=h(F) = \#\Pic^0(F)$ its class
number. We want upper and lower bounds for $h(F)$. Using the
Hasse-Weil theorem, one easily obtains $|\sqrt{q}-1|^{2g} \leq h
\leq |\sqrt{q}+1|^{2g}$. Unfortunately, the lower bound is only
useful when $q \geq 5$, and so for general $q$ we have the
following bound, which was shown to me by Henning Stichtenoth.

\begin{proposition}\label{Boundh}
We have
\[
h(F) \geq \frac{(q-1)(q^{2g} - 2gq^g + 1)}{2g(q^{g+1} - 1)}.
\]
\end{proposition}

\paragraph{Proof.}
We consider the constant field extension $F'=\FF_{q^{2g}} F$ of
$F$ of degree $2g$. The exact field of constants of $F'$ is
$\FF_{q^{2g}}$. Let $N'$ denote the number of rational (that is,
$\FF_{q^{2g}}$-rational) places of $F'$. The Hasse-Weil bound
gives us $N' \geq q^{2g} -2gq^g + 1.$ Let $Q$ be one such rational
place of $F'$, lying over the place $P$ of $F$. As $Q$ has degree
one, we get $2g = f(Q|P)\deg(P)$, and so $\deg(P)$ divides $2g$.
It follows that $(2g/\deg(P))\cdot P$ is an effective divisor of
degree $2g$ of $F$. As there are at most $2g$ places $Q$ above
$P$, we see that in this way we have constructed at least $N'/2g$
effective divisors of degree $2g$ of $F$. On the other hand, there
are exactly $h(q^{g+1}-1)/(q-1)$ such places, so we get
\[
\frac{h}{q-1}(q^{g+1}-1) \geq \frac{N'}{2g} \geq \frac{q^{2g}
-2gq^g + 1}{2g},
\]
from which the result follows.\Endproof

We now let $F=K=k(\sqrt{D})$ be an imaginary quadratic extension
of $k=\Fq(T)$, where $D\in A$ is square-free. Then the genus of
$K$ is given by
\[
g = \left\{\begin{array}{cl}
 (\deg(D) - 1)/2 & \mbox{if $\deg(D)$ is odd} \\
 (\deg(D) - 2)/2 & \mbox{if $\deg(D)$ is even.}\end{array}\right.
\]

Let $\OO$ be an order of conductor $f$ in $K$. Then, as in the
classical case, one may express $\#\Pic(\OO)$ in terms of $h(K)$
and $f$ (e.g. \cite{Rosen02}, Proposition 17.9), which, combined
with our bounds on $h(K)$, gives us
\begin{equation}\label{bound3}\label{bound4}
B_{\varepsilon}|Df^2|^{\frac{1}{2} - \varepsilon} \leq \#\Pic(\OO)
\leq C_{\varepsilon} |Df^2|^{\frac{1}{2} + \varepsilon}
\end{equation}
for every $\varepsilon>0$ and effectively computable positive
constants $B_{\varepsilon}$ and $C_{\varepsilon}$.

\subsection{Estimating the $j$-invariant}\label{sectionjestimate}

In this section we estimate the $j$-invariant using analytic
methods, following the first part of \cite{Brown92}. We point out
that later parts of that paper (the part concerning supersingular
reduction) have been shown to contain errors, but we will only use
results from the first (and supposedly correct) part.

\begin{definition}{\em
Let $z\in\Omega$. Then we define
\begin{eqnarray*}
|z|_A & = & \inf_{a\in A}|z - a|,\quad\mbox{and}\\
|z|_i & = & \inf_{x \in \kinf}|z - x|.
\end{eqnarray*}
The {\em imaginary modulus} $|z|_i$ plays the role of $|\Im(z)|$
in the classical case. }\end{definition}

Let $\phi$ be a CM Drinfeld module. Then
$\End(\phi)=\OO=A[\sqrt{d}]=A[f\sqrt{D}]$ is an order of conductor
$f$ in $K=k(\sqrt{D})$, where $D$ is the square-free part of $d =
Df^2$.

A non-zero ideal $\fa$ in $\OO$ is a rank 2 lattice in $\bC$. It follows
that $\fa$ is homothetic to the lattice $\Lambda_z = \langle z,1
\rangle$, for some $z\in\Omega$. This $z$ is determined up to
$\PGLA$-action, so we would like to have a fundamental domain for
this action. Unfortunately, a perfect analogue of the classical
fundamental domain for the $\SL(\ZZ)$-action on $\HH$ does not
seem to exist, but if we're only interested in $z\in\Omega$
quadratic over $k$, then we do have the next best thing.

\begin{definition}{\em
The {\em quadratic fundamental domain} is
\begin{eqnarray*}
\cD = \{ z\in\Omega  & | & \mbox{$z$ satisfies an equation of the form $az^2+bz+c=0$,}\\
 & & \mbox{where $a,b,c\in A$, $a$ is monic, $|b|<|a|\leq |c|$,}\\
 & & \mbox{and $\gcd(a,b,c)=1$}\}.
\end{eqnarray*}
}\end{definition}

In general we're only interested in $\cD\cap K$, which we denote
$\cD_K$. Then, as in the classical case, one may show that any
rank $2$ lattice in $K$ is homothetic to $\Lambda_z$ for some
$z\in\cD_K$. Moreover, we have

\begin{proposition}
If $z\in\cD_K$, then $|z|_i=|z|_A=|z|\geq 1$.
\end{proposition}

\paragraph{Proof.}
It suffices to show that $|z|_i=|z|\geq 1$. Write
$z=(-b+\sqrt{d})/2a$, where $d=b^2-4ac$. Then $|d|=|ac|\geq |a^2|$
and $|d|\leq |c^2|$. Hence $|z|=|\sqrt{d}/2a|\geq 1$.

We distinguish two cases.

{\bf (a)} If $\infty$ is ramified in $K/k$, then $\deg(d)$ is odd
and $v_{\infty}(\sqrt{d}/2a) \in \frac{1}{2}\ZZ \smallsetminus
\ZZ$ is half integral, so $|x| \neq |z|$ and $|z-x| \geq |z|
\;\;\forall x\in\kinf$.

{\bf (b)} If $\infty$ is inert, then $\deg(d)$ is even, but its
leading coefficient is not a square in $\Fq$. So
the leading coefficient of $\sqrt{d}/2a$ as a Laurent series in
$1/T$ is not in $\Fq$, and the leading terms of $x\in\kinf$ and
$\sqrt{d}/2a$ cannot cancel. Hence $|z-x| \geq |z|$ in this case,
too. \Endproof

Now we may estimate $|j(z)| = |j(\phi^z)|$.

\begin{theorem}\label{jestimateCM}
Suppose $q$ is odd. Let $z=(-b+\sqrt{d})/2a \in \cD_K$. Then
\begin{enumerate}
\item If $|z|=1$ then $|j(z)|\leq 1/q$. \item If $|z|>1$ then
$|j(z)|=B_q^{|z|},$ where
\[
B_q = \left\{\begin{array}{ll}q^q & \mbox{if $\deg(d)$ is even} \\
   q^{\sqrt{q}(q+1)/2} & \mbox{if $\deg(d)$ is odd.}\end{array}\right.
\]
\end{enumerate}
\end{theorem}

\paragraph{Proof.}
One just follows the proof of \cite{Brown92}, Theorem 2.8.2, using
the fact that $|z|_A=|z|_i=|z|\geq 1$ when $z\in\cD_K$. Then all
the calculations of \cite{Brown92} work and we do not need to
assume that $d$ be square-free (i.e. that $z$ correspond to a
Drinfeld module with complex multiplication by the full ring of
integers $\OO_K$ of $K$). \Endproof

\begin{corollary}\label{tojestimate}
Let $\OO$ be an order in $K$, and let $\fa\subset\OO$ be an
invertible ideal. Then $|j(\fa)|\leq  |j(\OO)|$, with equality if
and only if $\fa$ is principal.
\end{corollary}

\paragraph{Proof.}
Write $\OO=A[\sqrt{d}]$. Then the representatives
$\fa_1,\ldots,\fa_h$ of the ideal classes in $\Pic(\OO)$, with
$h=\#\Pic(\OO)$ and $\fa_1=\OO$, correspond to elements $z_i =
(-b_i+\sqrt{d_i})/2a_i\in\cD_K$, with $z_1=\sqrt{d_1}=\sqrt{d}$.
Now
\[
A[\sqrt{d_i}] = \End(\Lambda_{z_i}) = \End(j(\fa_i))= \OO =
A[\sqrt{d}]
\]
for every $i$, so we see that $d$ and $d_i$ differ only by the
square of a unit, hence can be assumed to be equal. Now we have
$|z_1|=|\sqrt{d}| > |\sqrt{d}/2a_i|=|z_i|$ for all $i\neq 1$ and
the result follows from Theorem \ref{jestimateCM}. (Note that if
$|a_i|=1$, then $a_i=1$ and $b_i=0$, so $i=1$.) \Endproof

It follows in particular that $j(\OO)$ is larger than any of its
other conjugates.

\subsection{CM heights}\label{sectionCMH}

\begin{definition}
{\em Let $\phi$ be a CM Drinfeld module, with
$\End(\phi)=A[\sqrt{d}]=A[f\sqrt{D}]$ and $j$-invariant
$j=j(\phi)$. Then we define the {\em CM height} of $\phi$ and of $j$ to be
\[
\HCM(j) = \HCM(\phi) = |d| = |f^2D|.
\]
If $x\in\bC$ is not a CM point, then we set $\HCM(x)=1$.
If $x=(x_1,\ldots,x_n)\in\AA^n(\bC)$ then we define
\[
\HCM(x) = \max\{\HCM(x_1),\ldots,\HCM(x_n)\}.
\]
}\end{definition}

This height is not to be confused with the usual notion of the
height of a Drinfeld module of finite characteristic. In fact, all
the Drinfeld modules here have generic characteristic. This
definition is analogous to the definition given in \cite{Breuer01}
for elliptic curves. The CM height is so-named because it forms a
true counting function on the CM points of $\AA^1(\bC)$ (and thus
also of $\AA^n(\bC)$).

\begin{proposition}
For every $\varepsilon >0$ we have
\[
\#\{ j\in\bC \;|\; \mbox{$j$ is CM and $\HCM(j)\leq t$}\} =
O(t^{3/2+\varepsilon}).
\]
\end{proposition}

\paragraph{Proof.}
For every order $\OO_d=A[\sqrt{d}]$, there are exactly
$\#\Pic(\OO_d)$ isomorphism classes of Drinfeld modules $\phi$
with $\End(\phi)=\OO_d$, namely those corresponding to the ideal
classes $[\fa]\in\Pic(\OO_d)$. So we have
\begin{eqnarray*}
  \#\{j\in\bC \;|\; \mbox{$j$ is CM and $\HCM(j)\leq t$}\}
 & = & \sum_{|d|\leq t}\#\Pic(\OO_d) \\
 & \leq & \sum_{|d| \leq t} C_{\varepsilon}|d|^{1/2 + \varepsilon} \qquad\mbox{from (\ref{bound4})} \\
 & = & O(t^{3/2 + \varepsilon}).
\end{eqnarray*}
\Endproof

Now that we may view the CM height as a height function, one may
ask how this compares to the usual (i.e. arithmetic) height in
$\PP^1$ (see \cite{Hindry-Silverman}, Part B, or \cite{LangFDG},
Chapter 3, for definitions). We have

\begin{proposition}\label{CMvsArith}
Let $j\in k^{sep}$ be a CM point, with $\End(j)=\OO=A[\sqrt{d}]$
and $\HCM(j)=|d|$. Then the (logarithmic) height of $j\in
\PP^1(k^{sep})$ is bounded by $h(j)\leq\HCM(j)^{1/2} + C_q$, where
\[
C_q = \left\{\begin{array}{ll}q & \mbox{if $\deg(d)$ is even} \\
  \sqrt{q}(q+1)/2 & \mbox{if $\deg(d)$ is odd.} \end{array}\right.
\]
\end{proposition}

\paragraph{Proof.}
Let $K = \OO\otimes_A k$ denote the CM field and set $F=K(j)$. We
recall that $j$ is integral over $A$, so that $|j|_v\leq 1$ for
any place $v$ of $F$ that does not lie over the (unique) place
$\infty$ of $K$. On the other hand, the place $\infty$ splits
completely in $F/K$, so for any place $v|\infty$ of $F$ we have
$|j|_v = |\sigma_v(j)|^2$, where $\sigma_v : F \hookrightarrow
\bC$ is the embedding of $F$ into $\bC$ corresponding to the place
$v$, and $|\cdot |$ denotes the unique absolute value of $\bC$.
This gives us
\begin{eqnarray*}
h(j) 
     & = & \frac{1}{2[F:K]}\sum_{v|\infty}\log_q(\max\{|j|_v,1\})\qquad\mbox{(as $j$ is integral)} \\
     & = & \frac{1}{2[F:K]}\sum_{\sigma\in\Gal(F/K)}\log_q(\max\{|\sigma(j)|^2,1\}) \\
     & \leq & \log_q|j(\OO)| \qquad\mbox{(from Corollary \ref{tojestimate})} \\
     & \leq & |z| + \log_q(B_q) \qquad\mbox{(from Theorem \ref{jestimateCM})} \\
     & = & |d|^{1/2} + C_q.
\end{eqnarray*}
The result follows. \Endproof

\subsection{CM points on curves}\label{AOcurvesSection}

We are now ready to prove our first main result: the effective
Andr\'e-Oort conjecture for the product of two Drinfeld modular
curves.

\paragraph{Proof of Theorem \ref{AOcurves}.}
We will prove Theorem \ref{AOcurves} for $n=2$, the extension to
general $n$ then follows by considering projections to pairs of
coordinates.

Let $X\subset\AA^2$ be a curve of degree $d$, as in Theorem
\ref{AOcurves}. Firstly, it is clear that the modular curves
$Y'_0(N)$ contain CM points of arbitrary height. We want to prove
the converse. Let $x=(x_1,x_2)\in X(\bC)$ be a CM point. From
Proposition \ref{CMvsArith} follows that it suffices to show that
$X$ is modular if $x$ has a large {\em CM height}. We may assume
that both projections $p_i : X \rightarrow \AA^1$ are dominant
(otherwise the result is trivial). We want to use Theorem
\ref{Topology}, so we must show that $X$ is stabilized by a
suitable Hecke operator.

Let $\OO_i = \End(x_i) = A[f_i\sqrt{D_i}]$ be orders of conductors
$f_i$ in the imaginary quadratic fields $K_i$, for $i=1,2$, and
let $K=K_1K_2$. Denote by $g_i$ the genus of $K_i$. Denote by
$F_s$ the separable closure of $k$ in $F$, 
and let $L$ be the Galois closure of $F_sK(x_1,x_2)$ over $k$.

Let $\fp$ be a prime of even degree in $k$ which splits completely
in $F_sK$ and does not divide $f_1f_2$. Let $\fP$ be a prime of
$L$ lying over $\fp$, and denote by $\fP_i$ its restriction to the
field $K_i(x_i)$.

From CM theory (Theorem \ref{CMtheory}) follows that
$\Gal(K_i(x_i)/K_i) \cong \Pic(\OO_i)$ and $\fp$ is unramified in
$L/k$. 
Denote by $\sigma\in\Aut(FL/FK)$ an extension of the Frobenius element
$(\fP,L/k)\in\Gal(L/k)$,
%
and let $\sigma_i = \sigma|_{K_i(x_i)} = (\fP_i,
K_i(x_i)/k) = (\fP_i, K_i(x_i)/K_i)$, as $\fp$ splits in $K_i$.
Moreover, we have cyclic isogenies $x_i \rightarrow
x_i^{\sigma_i}$ of degree $\fp$, so $(x_1,x_2)\in X\cap
T_{\AA^2,\fp}(X^{\sigma}) = X\cap T_{\AA^2,\fp}(X)$, as $\sigma$
acts trivially on $F$.

On the one hand, from Proposition \ref{degrees} follows that
$\deg(X\cap T_{\AA^2,\fp}(X)) \leq 4d^2(|\fp|+1)^2$. On the other
hand, 
the whole $\Gal(FK(x_1,x_2)/F)$-orbit of the point $(x_1,x_2)$ lies in
this intersection, and there are at least $\#\Pic(\OO_i)/m$ points
in this orbit (for $i=1$ and $i=2$). We must show that
$\#\Pic(\OO_i) > 4md^2(|\fp|+1)^2$, as then the intersection will
be improper, giving $X \subset T_{\AA^2,\fp}(X)$, as $X$ is
irreducible. Then the result will follow from Theorem
\ref{Topology}, if $|\fp|\geq\max(13,d)$.

It remains to show that there exists a prime $\fp$ which has the
desired properties. For this we use the \v{C}ebotarev Theorem (see
\cite{FriedJarden}, Proposition 5.16). Let $M$ be the 
Galois closure of $F_sK$ over $k$,
and set
\[
\pi_M(t)=\#\{ \fp \in A \;|\; \mbox{prime, split in $M$, and
$|\fp|=q^t$}\}.
\]
Let $\FF$ be the algebraic closure of $\Fq$ in $M$, let
$n_c=[\FF:\Fq]$ be the constant extension degree and $n_g = [M :
\FF k]$ be the geometric extension degree. If $n_c \nmid t$ then
$\pi_M(t)=0$. If $n_c | t$ then
\[
| \pi_M(t) - \frac{1}{n_g}q^t/t | < 4(g(M) + 2)q^{t/2}.
\]
Here $g(M)$ is the genus of $M$, which can be bounded with the
Castelnuovo inequality (\cite{Stichtenoth93}, Theorem III.10.3) 
to give $g(M)\leq C_1(g_1+g_2) +C_2g +C_3$,
and we also have $n_gn_c \leq C_4$. Here the $C_i$'s are
computable constants depending only on $m$.

Now we want both $\pi_M(t) > \deg(f_1f_2)=\log_q|f_1f_2|$ (so that
we have a split prime $\fp$ not dividing $f_1f_2$) and
$\#\Pic(\OO_i)
> 4md^2(q^t + 1)^2$ (so that $X\subset T_{\AA^2,\fp}(X)$).

Summarizing, we want a simultaneous solution $t\in 2n_c\NN$ to the
inequalities
\begin{equation}
\frac{1}{C_4}q^t/t - 4\big(C_1(g_1+g_2) + C_2g + C_3
+2\big)q^{t/2}
> \log_q|f_1f_2|
\end{equation}
\begin{equation}
B_{\varepsilon}(q^{g_i}|f_i|)^{1-\varepsilon} > 4md^2(q^t+1)^2
\quad\mbox{from (\ref{bound3})}
\end{equation}
for some $\varepsilon > 0$ and at least one of $i=1$ or $i=2$.

These inequalities hold with $q^t \geq \max(13,d)$ if
$\HCM(x_1,x_2) = \max(|D_1f_1^2|,|D_2f_2^2|)$ is larger than some
computable constant $B$, which depends on $d,m$ and $g$. \Endproof

As we may equip our Drinfeld modules with arbitrary level
structures - which play no role - we may replace each copy of
$\AA^1$ by a Drinfeld modular curve to obtain

\begin{corollary}\label{toAOcurves}
Let $X_1,\ldots,X_n$ be Drinfeld modular curves. Let
$Z=X_1\times\cdots\times X_n$, and let $X \subset Z$ be an
irreducible algebraic curve. Then the following are equivalent:
\begin{enumerate}
\item $X$ contains infinitely many CM points \item $X$ contains at
least one CM point of height larger than some effectively
computable constant which depends only on $Z$, $\deg(X)$ and the
field of definition of $X$. \item There exists a non-empty subset
$S\subset \{1,\ldots,n\}$ for which we may write
\[ Z\cong Z_S\times Z'_S =
   \Big(\prod_{i \in S}X_i\Big) \times \Big(\prod_{i \notin S}X_i\Big),
\]
\[ X=X' \times \Big(\prod_{i \notin S} \{x_i\}\Big), \]
where the $x_i\in X_i$ are CM points and $X'$ is a special curve
in $Z_S$.
\end{enumerate}
\end{corollary}

\subsection{CM points on varieties}\label{AOvarietiesSection}

In this last section we prove our other main result: the
Andr\'e-Oort conjecture for subvarieties of the product of $n$
Drinfeld modular curves.

\paragraph{Proof of Theorem \ref{AOvarieties}.}
As CM points are defined over $k^{sep}$, so is $X$. Hence there
exists a finite Galois extension $F$ of $k$ such that $X$ is
defined over $F$.

Set $d=\dim(X)$. We will use induction on $d$. From Theorem
\ref{AOcurves} and Corollary \ref{toAOcurves} we know that the
result already holds for $d=1$. We now suppose $d\geq 2$, $n\geq
3$ and that the result is already known for dimensions less than
$d$. We claim that we may assume that $X\subset \AA^n$ is a
hypersurface.

Indeed, $X$ is an irreducible component of
\[
\bigcap_{\stackrel{I \subset \{1,\ldots,n\}}{\#I = d+1}}
p_I^{-1}(p_I(X)).
\]
The CM points are Zariski-dense in the hypersurfaces $p_I(X)
\subset\AA^I$, and if these are modular, then so are the
$p_I^{-1}p_I(X)$, and thus also $X$. This proves our claim.
Furthermore, we may assume that all the projections $p_i : X
\rightarrow \AA^1$ are dominant.

For a given constant $B>0$ we may assume that every point
$x=(x_1,\ldots,x_n)\in S$ satisfies $\HCM(x_i)>B$ for all
$i=1,\ldots,n$, as the set
\[
\{x \in X(\bC) \;|\; \mbox{$\HCM(x_i) \leq B$ for some
$i=1,\ldots,n$}\}
\]
is contained in a proper closed subvariety of $X$.

\paragraph{Step 1.}
Choose a point $x=(x_1,\ldots,x_n)\in S$. Suppose that we have
primes $\fp_1,\ldots,\fp_{d-1}$ of $k$ of even degree, 
satisfying the following conditions:
\begin{enumerate}
\item Each $\fp_j$ splits completely in every $\OO_i=\End(x_i)$
for $i=1,\ldots,n$ and in $F$. \item $|\fp_1|\geq\max\{13,\deg
X\}$ \item $|\fp_{j+1}| \geq (\deg
X)^{2^j}\prod_{m=1}^j(2|\fp_m|+2)^{n2^{j-m}}$ for $j=1,\ldots,d-2$
\item We have $\#\Pic(\OO_i) > [F:k](\deg
X)^{2^{d-1}}\prod_{m=1}^{d-1}(2|\fp_m|+2)^{n2^{d-m-1}}$ for each
$i=1,\ldots,n$, for which it suffices to assume
\end{enumerate}
\begin{equation}\label{Veq1}
 \#\Pic(\OO_i) > [F:k]|\fp_{d-1}|^2(2|\fp_{d-1}|+2)^n.
\end{equation}
Then, as in the proof of Theorem \ref{AOcurves}, it follows that
\[
\Gal(F^{sep}/F)\cdot x \subset X\cap T_{\AA^n,\fp_j}(X),\quad
j=1,\ldots,d-1.
\]

Let $X_1$ be an $F$-irreducible component of $X\cap T_{\fp_1}(X)$
containing $x$. Now either $X_1 = X$, in which case $X_1 \subset
T_{\fp_1}(X_1)$ and $X_1$ is modular (Theorem \ref{Topology2}), or
$\dim(X_1) < \dim(X)$. In the latter case we repeat the procedure:
We let $X_2$ be an $F$-irreducible component of $X_1\cap
T_{\fp_2}(X_1)$ containing $x$, and so on. We thus produce a
sequence $X_1,X_2,\ldots$ of $F$-irreducible subvarieties of $X$
of strictly decreasing dimension. But, as the $X_j$ are defined
over $F$, the full $\Gal(F^{sep}/F)$-orbit of $x$ is contained in
each $X_j$. Moreover, after at most $d-1$ steps we arrive at
$\dim(X_j)\leq 1$, and
\begin{eqnarray*}
\deg X_j  & \leq & (\deg X)^{2^j}
\prod_{m=1}^j(2|\fp_m|+2)^{n2^{j-m}}
                   \quad\mbox{(using Proposition \ref{degrees})} \\
          & < &  \#\Pic(\OO_i)/[F:k] \leq \#\Gal(F^{sep}/F)\cdot x
                   \quad\mbox{(as $j\leq d-1$).}
\end{eqnarray*}
Hence $X_j$ must have dimension at least $1$. In summary, this
process must terminate, after at most $d-1$ steps, with some $X_j$
of dimension at least $1$, satisfying $X_j \subset
T_{\fp_{j+1}}(X_j)$. Hence $X_j$ is modular.

By varying $x\in S$, we see that we have covered $X$ by a
Zariski-dense family of modular subvarieties $X_x$ for $x\in S$.
We now show that the $X_x$'s are in fact pure modular. Suppose
not. Recall that each $X_x$ is $F$-irreducible, so if it's not
pure then it is the $\Gal(F^{sep}/F)$-orbit of a modular variety
of the form $Y_x \times \{y_x\}$, where $y_x$ is a CM point (in
fact a projection of $x$) and $Y_x$ a pure modular variety. But as
the $\Gal(F^{sep}/F)$-orbit of the point $y_x$ is larger than the
degree of $X_x$, by construction, this would mean that the number
of (geometrically) irreducible components of $X_x$ of maximal
dimension is larger than $\deg(X_x)$, which is impossible. So each
$X_x$ is in fact pure modular. Now each $X_x$ contains a
Zariski-dense family of pure modular curves, hence so does $X$.

\paragraph{Step 2.}
We want to show that $X$ is modular, using the fact that $X$
contains a Zariski-dense family of pure modular curves. For ease
of notation we will denote this family by $S$ and the pure modular
curves by $s\in S$.

Choose a CM point $x_1 \in \AA^1(\bC)$ and consider the
intersection
\[
X_1 = X \cap (\{x_1\}\times\AA^{n-1}).
\]
As each curve $s\in S$ is pure modular, it intersects $X_1$ in at
least one CM point. We denote by $X'$ the Zariski closure of these
points:
\[
X' = \overline{\cup_{s\in S}(s\cap X_1)}^{Zar}.
\]
Now if $\dim(X')=\dim(X)$, then $X\subset\{x_1\}\times\AA^{n-1}$,
which is impossible: we had assumed in the beginning that all
projections $p_i : X \rightarrow \AA^1$ are dominant. Hence
$\dim(X') < \dim(X)$. Then it follows from the induction
hypothesis that all the (geometrically) irreducible components of
$X'$ are modular. Write $X'=X_1'\cup \cdots \cup X_r'$ as the
union of $r$ irreducible components. Then the points of $s\cap
X_1$ distribute amongst these components. By restricting $S$ to a
Zariski-dense subfamily, and renumbering the components of $X'$,
we may assume that $X_1'$ contains at least $1/r$ of the points of
$s\cap X_1$ for every $s\in S$.

If, up to permutation of coordinates, $X_1'$ is of the form
$\{y\}\times\AA^m$ for some $m<n-1$ and $y$ a CM point in
$\AA^{n-m}$, then it follows that $X$ is of the form (again up
to permutation of coordinates) $Y \times \AA^{n-2}$, where $Y$ is
an irreducible curve in $\AA^2$. But then $Y$ contains infinitely
many CM points, hence is modular. In this case we see that $X$ is
modular.

So we may now assume that at least one modular curve appears as a
factor of $X_1'$. Then there exists some pair of coordinates
$1<i<j$ such that
\[
p_{i,j}(X_1') = Y'_0(m),
\]
for some {\em fixed} $m\in A$ .

Let $s'=p_{\{1,i,j\}}(s)\subset\AA^3$ be characterized by the
triplet $(N_{s,1},N_{s,i},N_{s,j}) \in A^3$ as in
\S\ref{CurveSection}, and assume, by restricting $S$ to a
Zariski-dense subfamily and permuting coordinates, that we always
have $|N_{s,i}| \leq |N_{s,j}|$. Fix $s\in S$ and fix also $x_i$
such that we have a point $(x_1,x_i,x_j) \in s'$. We want to find
many points $x_j$ with this property.

For each prime $\fp\in A$, consider the tree $\TT_{\fp}$. Then a
{\em generic} point $(x_1,x_i,x_j)$ of $s'$ corresponds to a
triple of vertices $(v_{\fp,1},v_{\fp,i},v_{\fp,j})$, at distances
$(n_{\fp,1},n_{\fp,i},n_{\fp,j})$ from the center $v_{\fp,c}$. The
family of vertices $(v_{\fp,c})_{\fp\in A}$ corresponds to a point
$x_c\in\bC$ which we call the {\em center} of $(x_1,x_i,x_j)$.
The possible choices of $v_{\fp,j}$ correspond to the length
$n_{\fp,j}$ paths leading out from $v_{\fp,c}$ and disjoint from
the two paths leading to $v_{\fp,1}$ and $v_{\fp,i}$. This gives
$(|\fp|-1)(|\fp|+1)^{n_{\fp,j}-1}$ possibilities if $n_{\fp,j}\geq
1$ (and just one if $n_{\fp,j}=0$). Multiplying over all primes
$\fp\in A$ then shows that there are $\prod_{\fp|
N_{s,j}}(|\fp|-1)(|\fp|+1)^{n_{\fp,j}-1}$ possible choices for
$x_j$. So we have counted the number of suitable cyclic degree
$N_{s,j}$ isogenies from $x_c$ to $x_j$. But $(x_1,x_i,x_j)$ is a
CM point, not a generic point, so some of these isogenies will
produce the same point $x_j$, corresponding to non-trivial
endomorphisms
$\alpha\in \End(x_c)$ of norm $N_{K/k}(\alpha)=N_{s,j}^2$. The
number of such endomorphisms is at most
$\prod_{\fp|N_{s,j}}(2n_{\fp,j}+1)$, so the number of distinct
values of $x_j$ satisfying $(x_1,x_i,x_j)\in s'$ tends to infinity
as $N_{s,j}$ increases.

But $1/r$ of these points also satisfy $(x_i,x_j)\in Y'_0(m)$, of
which there can be at most $\psi(m)$, for fixed $x_i$. So we have
shown that $N_{s,j}$, and thus also $N_{s,i}$, is bounded as $s$
ranges through $S$.

It follows that there are only finitely many possibilities for
$p_{i,j}(s)=Y'_0(N_{s,i}N_{s,j})$. By replacing $S$ with a
Zariski-dense subfamily, we may assume there is only one:
$p_{i,j}(s) = Y'_0(N_0)$ for all $s\in S$. Now, after a
permutation $(i,j)\mapsto (n-1,n)$ of coordinates, we see that
\begin{eqnarray*}
S & \subset & \AA^{n-2}\times Y'_0(N_0), \quad\mbox{and so} \\
X = \overline{S}^{Zar} & \subset & \AA^{n-2}\times Y'_0(N_0).
\end{eqnarray*}
But $X$ is a hypersurface, so we have in fact $X = \AA^{n-2}\times
Y'_0(N_0)$, which is modular. This is what we set out to prove.

\paragraph{Step 3.}
It remains to show that we can find primes $\fp_j$ with the
desired properties. Recall that $x=(x_1,\ldots,x_n)$ and each
$\OO_i=\End(x_i)$ is an order of conductor $f_i$ in the imaginary
quadratic field $K_i$ of genus $g_i$.

Set $|\fp_j|=q^{t_j}$ for $j=1,\ldots,d-1$. Firstly, we need
(\ref{Veq1}), which, combined with the lower bound for the class
number (\ref{bound3}), gives
\begin{equation}\label{Veq2}
B_{\varepsilon}(q^{g_i}|f_i|)^{1-\varepsilon} >
[F:k]q^{2t_{d-1}}(2q^{t_{d-1}}+2)^n.
\end{equation}
Secondly, the $\fp_j$'s must be well spaced out, i.e. we need
\begin{equation}\label{Veq3}
q^{t_{j+1}} \geq (\deg
X)^{2^j}\prod_{m=1}^j(2q^{t_m}+2)^{n2^{j-m}}.
\end{equation}
Thirdly, each $\fp_j$ must split completely in $FK$, where
$K=K_1\cdots K_n$, and not divide $f_1\cdots f_n$. Here the
\v{C}ebotarev Theorem says
\[
|\pi_{FK}(t_j) - \frac{1}{n_g}q^{t_j}/t_j | <
4(g(FK)+2)q^{t_j/2}\quad\mbox{and $n_c|t_j$},
\]
where $n_g$ denotes the geometric extension degree of $FK/k$,
$n_c$ denotes the constant extension degree, and $g(FK)$ is the
genus of $FK$. We have $n_gn_c = [FK:k] \leq 2^n[F:k]$ and we may
bound $g(FK)$ from above via the Castelnuovo Inequality to obtain
$g(FK) \leq C_1(g_1+\cdots+g_n) + C_2$ for some computable
constants $C_1$ and $C_2$ depending on the field $F$.
We need $\pi_{FK}(t_j) > \log_q|f_1\cdots f_n|$, to obtain a split
prime $\fp_j$ that does not divide any of the conductors.

In summary, we want $d-1$ simultaneous solutions
$t_1,\ldots,t_{d-1}\in 2n_c\NN$ to the inequalities $q^{t_1} \geq \max(13, \deg
X)$, (\ref{Veq2}), (\ref{Veq3}) and
\begin{equation}
 \frac{1}{2^n[F:k]} q^{t_j}/t_j - 4\big(C_1(g_1+\cdots+g_n)+C_2+2\big)q^{t_j/2}
 > \log_q|f_1\cdots f_n|. \label{Veq4}
\end{equation}

If we choose the constant $B$ sufficiently large then, as $B <
\HCM(x_i) \leq q^{2g_i+1}|f_i|^2$ for all $i=1,\ldots,n$, such a
set of solutions $(t_1,\ldots,t_{d-1})$ exists. \Endproof

\begin{corollary}\label{toAOvarieties}
Let $X_1,\ldots,X_n$ be Drinfeld modular curves. Let
$Z=X_1\times\cdots\times X_n$, and let $X \subset Z$ be an
irreducible algebraic subvariety. Then the following are
equivalent:
\begin{enumerate}
\item $X$ contains a Zariski-dense set of CM points \item There
exists a partition $\{1,\ldots,n\} = \coprod_{i=0}^{g+1}S_i$ for
which we may write
\[
 Z\cong \prod_{i=0}^{g+1} Z_i = \prod_{i=0}^{g+1}\Big(\prod_{j \in S_i}X_j\Big),
\]
\[
X = Z_0 \times \prod_{i=1}^g Y_i \times \{x\},
\]
where each $Y_i$ is a special curve in $Z_i$ (for $i=1,\ldots,g$)
and $x$ is a CM point in $Z_{g+1}$.
\end{enumerate}
\end{corollary}

We remark that Corollary \ref{toAOvarieties} with $X_1 = \cdots = X_n = X_0(MN)$
has an application to Heegner points on elliptic curves over $k$ with
conductor $N\cdot\infty$, see \cite{Breuer04}.

\begin{center}
\rule[7mm]{50mm}{0.2pt}

Department of Mathematics,\\
University of Stellenbosch,\\
Stellenbosch 7600, South Africa.\\
{\em email:} fbreuer@sun.ac.za



\end{center}

\end{document}